\documentclass{article}

\usepackage{arxiv}
\usepackage{graphicx, color}
\usepackage{amsmath,amssymb}
\usepackage[utf8]{inputenc} % allow utf-8 input
\usepackage[T1]{fontenc}    % use 8-bit T1 fonts
\usepackage{hyperref}       % hyperlinks
\usepackage{url}            % simple URL typesetting
\usepackage{booktabs}       % professional-quality tables
\usepackage{amsfonts}       % blackboard math symbols
\usepackage{microtype}      % microtypography
\usepackage{lipsum}

\title{Addition formulas of Leaf Functions and Hyperbolic Leaf Functions}

\author{
  Kazunori Shinohara\thanks{10-3 Takiharu-cho, Minami-ku, Nagoya 457-8530, Japan} \\
  Department of Mechanical Systems Engineering\\
  Daido University\\
  10-3 Takiharu-cho, Minami-ku, Nagoya 457-8530, Japan \\
  \texttt{shinohara@06.alumni.u-tokyo.ac.jp} \\
}

\begin{document}
\maketitle

\begin{abstract}
Addition formulas exist in trigonometric functions. Double-angle and half-angle formulas can be derived from these formulas. Moreover, the relation equation between the trigonometric function and the hyperbolic function can be derived using an imaginary number. The inverse hyperbolic function $\mathrm{arsinh}(r)=\int_{0}^{r} \frac{1}{\sqrt{1+t^2} }\mathrm{d}t$ is similar to the inverse trigonometric function $\mathrm{arcsin}(r)=\int_{0}^{r} \frac{1}{\sqrt{1-t^2} }\mathrm{d}t$, such as  the second degree  of a polynomial and the constant term 1, except for the sign $-$ and $+$. Such an analogy holds not only when the degree of the polynomial is 2, but also for higher degrees. As such, a function exists with respect to the leaf function through the imaginary number $i$, such that the hyperbolic function exists with respect to the trigonometric function through this imaginary number. In this study, we refer to this function as the hyperbolic leaf function. By making such a definition, the relation equation between the leaf function and the hyperbolic leaf function makes it possible to easily derive various formulas, such as addition formulas of hyperbolic leaf functions based on the addition formulas of leaf functions. Using the addition formulas, we can also derive the double angle and half-angle formulas. We then verify the consistency of these formulas by constructing graphs and numerical data. 
\end{abstract}

\section{Introduction}
\label{Introduction}
\subsection{Leaf Functions and Hyperbolic Leaf Functions}
\label{Leaf Functions and Hyperbolic Leaf Functions}
An ordinary differential equation consists of both a function raised to the $2n-1$ power and the second derivative of the function.

\begin{equation}
\frac{\mathrm{d}^2r(l)}{\mathrm{d}l^2}=-nr(l)^{2n-1} \label{1.1.1}
\end{equation}

The preceding equation is the ODE that motivated this study. 
Although the equation (\ref{1.1.1}) is a simple ordinary differential equation, it has a very important meaning because it generates characteristic waves. By numerically analyzing the solution that satisfies this equation, we can obtain regular and periodic waves\cite{Kaz_cl, Kaz_sl}. The form of these waves differs from the form of the waves based on trigonometric functions. The function that satisfies this ordinary differential equation is called a leaf function, and it describes the features of these functions. Eq. (\ref{1.1.1}) is transformed as follows:
\newpage

\begin{equation}
l=\int_{0}^{r} \frac{\mathrm{d}t}{\sqrt{1-t^{2n}}}(=\mathrm{arcsleaf}_n(r)) \label{1.1.2}
\end{equation}

The preceding integral is defined as the inverse function $\mathrm{arcsleaf}_n (l)$ of the leaf function. Another function can be defined as follows:

\begin{equation}
l=\int_{r}^{1} \frac{\mathrm{d}t}{\sqrt{1-t^{2n}}}(=\mathrm{arccleaf}_n(r)) \label{1.1.3}
\end{equation}

The preceding integral is also defined as the inverse function $\mathrm{arccleaf}_n(r)$ of the leaf function with a different integral domain compared to Eq. (\ref{1.1.2}). The variable $n$ represents a natural number, and it is referred to as the basis. Moreover, the ordinary differential equation that is satisfied by the hyperbolic functions $r(l) = \mathrm{sinh}(l)$ and $r(l) = \mathrm{cosh}(l)$ is described as follows.

\begin{equation}
\frac{\mathrm{d}^2r(l)}{\mathrm{d}l^2}=r(l) \label{1.1.4}
\end{equation}

Compared to Eq. (\ref{1.1.1}), the difference in Eq. (\ref{1.1.4}) is the positive sign on the right hand side of the equation. The inverse hyperbolic functions $\mathrm{arsinh}(r)$ and $\mathrm{arcosh}(r)$ are well known as:

\begin{equation}
l=\int_{0}^{r} \frac{\mathrm{d}t}{\sqrt{1+t^{2}}}(=\mathrm{asinh}(r)) \label{1.1.5}
\end{equation}

\begin{equation}
l=\int_{1}^{r} \frac{\mathrm{d}t}{\sqrt{t^{2}-1}}(=\mathrm{acosh}(r)) \label{1.1.6}
\end{equation}

The contents of the root in the integrand constitute a polynomial. The polynomial of the inverse hyperbolic function and that of the inverse trigonometric function both have a degree of 2. The magnitude $1$ of the constant term in the root is also the same. The difference between the inverse functions of the trigonometric function and the hyperbolic function is the sign ($''+''$ and $''-''$) of the polynomial in the root. Using Eqs. (\ref{1.1.5}) and (\ref{1.1.6}), it is seen that trigonometric functions and hyperbolic functions have relational equation through imaginary numbers. Based on this relationship, similar functions also could be paired with leaf functions though analogy relation (see Appendix D in detail). These functions are called hyperbolic leaf functions and consist of two functions. One function is defined as follows.

\begin{equation}
r(l)=\mathrm{sleafh}_n(l)(n=1,2,3 \cdots) \label{1.1.7}
\end{equation}

The limit exists for the function $\mathrm{sleafh}_n(l)$ (see Appendix F). The domain of the variable $l$ is defined as follows:

\begin{equation}
-\zeta_n<l<\zeta_n \label{1.1.8}
\end{equation}

The initial conditions of the preceding equation are defined as follows.

\begin{equation}
r(0)=\mathrm{sleafh}_n(0)=0 (n=1,2,3 \cdots) \label{1.1.9}
\end{equation}

\begin{equation}
\frac{\mathrm{d}r(0)}{\mathrm{d}l}=\frac{\mathrm{d}}{\mathrm{d}l} \mathrm{sleafh}_n(0)=1 (n=1,2,3 \cdots) \label{1.1.10}
\end{equation}

Next, the another function is defined as follows:

\begin{equation}
r(l)=\mathrm{cleafh}_n(l) (n=1,2,3 \cdots) \label{1.1.11}
\end{equation}

The limit exists for the function $\mathrm{cleafh}_n(l)$ (see Appendix G). The domain of the variable $l$ is as follows:

\begin{equation}
-\eta_n<l<\eta_n \label{1.1.12}
\end{equation}

The initial conditions of the preceding equation are defined as follows.

\begin{equation}
r(0)=\mathrm{cleafh}_n(0)=1 (n=1,2,3 \cdots) \label{1.1.13}
\end{equation}

\begin{equation}
\frac{\mathrm{d}r(0)}{\mathrm{d}l}=\frac{\mathrm{d}}{\mathrm{d}l} \mathrm{cleafh}_n(0)=0 (n=1,2,3 \cdots) \label{1.1.14}
\end{equation}

The ordinary differential equations that are satisfied by the hyperbolic leaf functions that correspond to both equation (\ref{1.1.7}) and equation (\ref{1.1.11}) are as follows.

\begin{equation}
\frac{\mathrm{d}^2r(l)}{\mathrm{d}l^2}=nr(l)^{2n-1} \label{1.1.15}
\end{equation}

The inverse function of the hyperbolic leaf function is derived as follows:

\begin{equation}
l=\int_{0}^{r} \frac{\mathrm{d}t}{\sqrt{1+t^{2n}}}(=\mathrm{asleafh}_n(r)) (n=1,2,3 \cdots) \label{1.1.16}
\end{equation}

\begin{equation}
l=\int_{1}^{r} \frac{\mathrm{d}t}{\sqrt{t^{2n}-1}}(=\mathrm{acleafh}_n(r)) (n=1,2,3 \cdots) \label{1.1.17}
\end{equation}

Here, the prefix $a$ of both hyperbolic leaf functions $\mathrm{sleafh}_n (l)$ and $\mathrm{cleafh}_n(l)$ are defined as the inverse functions.

\subsection{Comparison of Legacy functions}
\label{Comparison}
The leaf functions and the hyperbolic leaf functions based on the basis $n=1$ are as follows:
\begin{equation}
\mathrm{sleaf}_{1}(l)=\mathrm{sin}(l) \label{1.1.18}
\end{equation}
\begin{equation}
\mathrm{cleaf}_{1}(l)=\mathrm{cos}(l) \label{1.1.19}
\end{equation}
\begin{equation}
\mathrm{sleafh}_{1}(l)=\mathrm{sinh}(l) \label{1.1.20}
\end{equation}
\begin{equation}
\mathrm{cleafh}_{1}(l)=\mathrm{cosh}(l) \label{1.1.21}
\end{equation}

Lemniscate functions were proposed by Johann Carl Friedrich Gauss\cite{Gauss}. The relation equations between these functions and leaf functions are as follows:

\begin{equation}
\mathrm{sleaf}_{2}(l)=\mathrm{sl}(l) \label{1.1.22}
\end{equation}
\begin{equation}
\mathrm{cleaf}_{2}(l)=\mathrm{cl}(l) \label{1.1.23}
\end{equation}
\begin{equation}
\mathrm{sleafh}_{2}(l)=\mathrm{slh}(l) \label{1.1.24}
\end{equation}

The definition of the function $\mathrm{slh}(l)$ in Eq.(\ref{1.1.24}) can be confirmed based on references\cite{Car, Neu}. A function corresponding to the hyperbolic leaf function $\mathrm{cleafh}_2(l)$ is not described in the literature\cite{Wei}. In the case where the basis $n \geqq 3$, the leaf function or the hyperbolic leaf function cannot be represented by a legacy function such as the lemniscate function.

\subsection{Originality and Purpose}
\label{Originality and Purpose}
The purpose of this report is to propose addition formulas for the hyperbolic leaf functions with basis $n=2$ and $n=3$, in addition to establishing both double-angle and the half-angle formulas using addition formulas. A similar analogy exists in the relation between the leaf function and the hyperbolic leaf function such that the relation between the trigonometric function and the hyperbolic function can be derived using imaginary numbers. Using this analogy, the addition formulas of hyperbolic leaf functions based on $n = 3$ can be derived from the addition formulas of leaf functions based on $n = 3$. Using addition formulas, we present numerical data and curves derived from the hyperbolic leaf function and show that these addition formulas in the section 2 are consistent.

\subsection{Contribution}
\label{sec:Contribution}
 The leaf functions are closely related to the Jacobi elliptic functions. The Jacobi elliptic function originated from the lemniscate function. In 1691, Jacob Bernoulli noticed that the arc length OP of the lemniscate curve was the same as the integral of the equation\cite{Roy}. As shown in Fig. \ref{fig1.1}, $l$ represents the length of arc OP. The arc OP is represented as:
 
\begin{figure*}[tb]
\begin{center} 
\includegraphics[width=0.75 \textwidth]{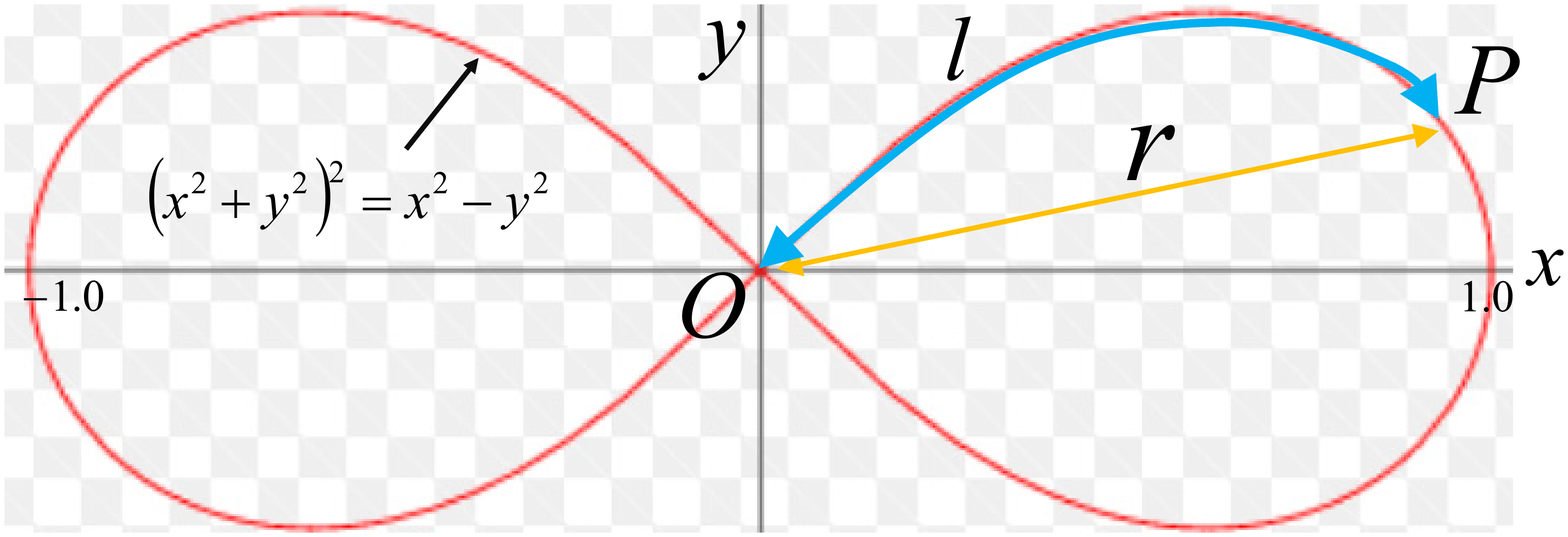}
\caption{ Lemniscate of Bernoulli }
\label{fig1.1}
\end{center}
\end{figure*}

\begin{equation}
arc OP=\int_{0}^{r} \frac{\mathrm{d}t}{\sqrt{1-t^4}}(=l) \label{1.1.25}
\end{equation}

The curve in Fig. \ref{fig1.1} can be expressed using variables $x$ and $y$ as:

\begin{equation}
(x^2+y^2)^2=x^2- y^2 \label{1.1.26}
\end{equation}

When the basis of the leaf function, n = 2, the curve for the leaf function is the same as the lemniscate curve. In 1718, Giulio Carlo de' Toschi di Fagnano published a paper explaining how arc OP could be divided into two equal parts using only a straightedge and a compass\cite{Ayo}. He discovered that the length $r$ was twice the length $u$, as shown in  Fig. \ref{fig1.2}. This led to the derivation of the addition theorem for the lemniscate function:

\begin{figure*}[tb]
\begin{center} 
\includegraphics[width=0.75 \textwidth]{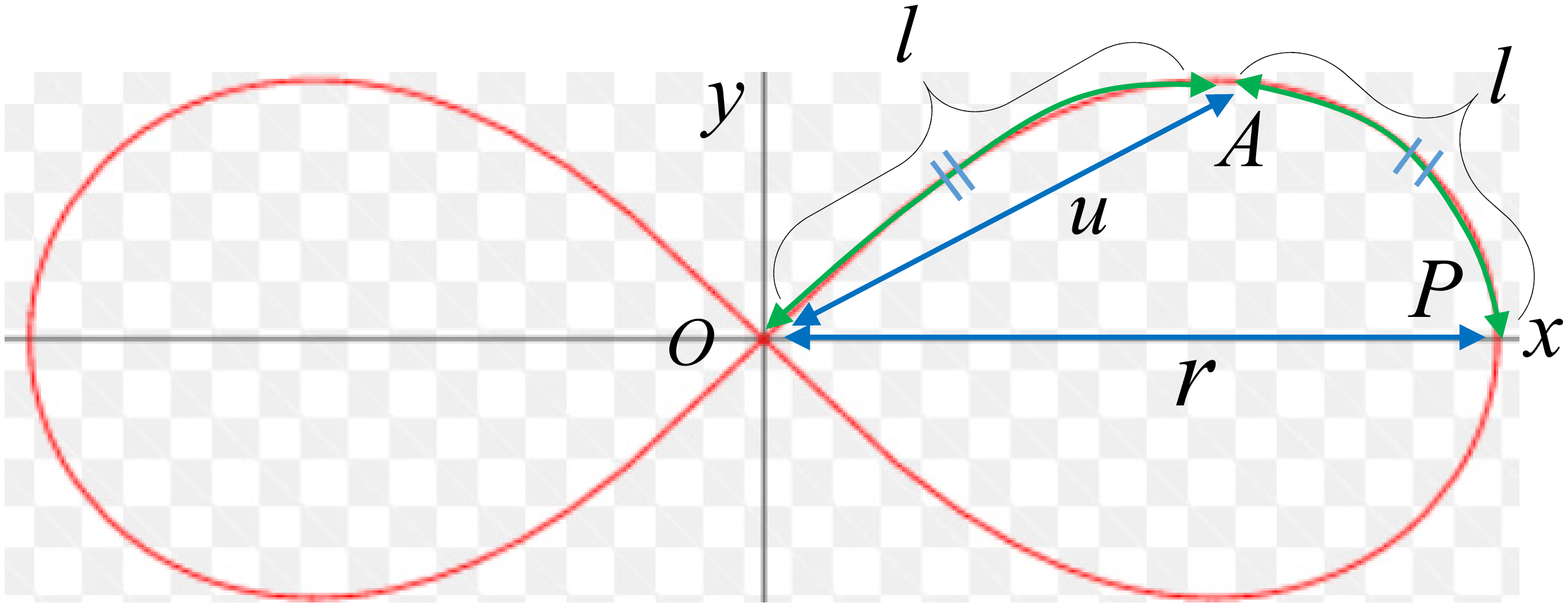}
\caption{ Fagnano and the Lemniscate }
\label{fig1.2}
\end{center}
\end{figure*}

\begin{equation}
r=\frac{ 2 u \sqrt{1- u^4}}{1+u^4} \: \mathrm{under} \: \mathrm{the} \: \mathrm{condition} \:  arc OP = 2 arc OA \label{1.1.27}
\end{equation}

\begin{equation}
arc OA= \int_{0}^{u} \frac{\mathrm{d}t}{\sqrt{1-t^4}} \label{1.1.28}
\end{equation}

After reading Fagnano's paper, Leonhard Euler found the addition formula for the lemniscate function in 1752\cite{Euler}. In the formula, the sum of the integral forms of arbitrary variables $u$ and $v$ equals the integral form $r$:

\begin{equation}
r=\frac{ u \sqrt{1-v^4} + v \sqrt{1-u^4}  }{ 1+ u^2 v^2 } \label{1.1.29}
\end{equation} 

The above-mentioned relation satisfies the integral equation as follows:

\begin{equation}
\int_{0}^{r} \frac{\mathrm{d}t}{\sqrt{1-t^4}} = \int_{0}^{u} \frac{\mathrm{d}t}{\sqrt{1-t^4}} + \int_{0}^{v} \frac{\mathrm{d}t}{\sqrt{1-t^4}}  \label{1.1.30}
\end{equation}

In 1796, Carl Friedrich Gauss derived the addition formula using the lemniscate function\cite{Gauss}. The inverse lemniscate function is expressed as:

\begin{equation}
\mathrm{arcsl}(r)=\int_{0}^{r} \frac{\mathrm{d}t}{\sqrt{1-t^4}} (=l_1 \: or \: l_2) \label{1.1.31}
\end{equation} 

Using this definition of inverse arc sl, the addition formula of lemniscate function $\mathrm{sl}(l_1+l_2)$ is derived as:

\begin{equation}
\mathrm{sl}(l_1+l_2)=\frac{ sl(l_1) \sqrt{1-sl(l_2)^4} + sl(l_2) \sqrt{1-sl(l_1)^4}  }{ 1+ sl(l_1)^2 sl(l_2)^2 } \label{1.1.32}
\end{equation} 

Phases $l_1$ or $l_2$ of lemniscate function sl can be extended to complex variables $i \cdot l_1$ or $i \cdot l_2$, respectively. In 1827, Carl Gustav Jacob Jacobi derived the inverses of the Jacobi elliptic functions\cite{Jacobi}. To derive the formula, the term $t^2$ is added to the root of the integrand denominator in Eq. (\ref{1.1.31}):

\begin{equation}
\mathrm{arcsn}(r,k)=\int_{0}^{r} \frac{\mathrm{d}t}{\sqrt{1-(1+k^2)t^2+k^2 t^4}} \label{1.1.33}
\end{equation} 

Eq. (\ref{1.1.33}) represents the inverse Jacobi elliptic function sn, where $k$ is a constant. There exist 12 Jacobi elliptic functions, including cn and dn. In the Eq. (\ref{1.1.33}), the variable $t$ is to the fourth power in the denominator. Jacobi did not discuss the variable $t$ to higher powers, such as follows:

\begin{equation}
\int_{0}^{r} \frac{\mathrm{d}t}{\sqrt{1-t^6}} \: , \: \int_{0}^{r} \frac{\mathrm{d}t}{\sqrt{1-t^8}}\: , \: \int_{0}^{r} \frac{\mathrm{d}t}{\sqrt{1-t^{10} }} \ldots\  \label{1.1.34}
\end{equation}

In other words, there had been no discussion for $n = 3$ in Eqs. (\ref{1.1.2}), (\ref{1.1.3}), (\ref{1.1.16}), and (\ref{1.1.17}). Therefore, the addition formulas for the leaf function were investigated for n=3\cite{Kaz_add}. In case of $n\geqq3$, no clear description for the addition formulas of hyperbolic leaf functions exists. On the contrary, $n=1$ represents hyperbolic functions sinh($l$) and cosh($l$). Therefore, the addition formulas of the hyperbolic leaf function and hyperbolic function are the same. The hyperbolic leaf function with $n=2$ represents hyperbolic lemniscate function slh($l$). No clear description of the addition formulas of function slh($l$) exists.

\subsection{Advantage and Disadvantage}
\label{Advantage and Disadvantage}
In physics, the nonlinear duffing equation represents a model for the spring pendulum whose spring stiffness does not obey Hooke's law. This undamped duffing equation is represented as:

\begin{equation}
\mathrm{Cubic-Quintic \: Duffing \: Equation:} \frac{\mathrm{d}^2r}{\mathrm{d}l^2}+ \alpha r+ \beta r^3 + \mu r^5=0 \label{1.1.35}
\end{equation}

To solve the above equation, numerical analysis or analytical approximate solutions have been applied\cite{Sibanda, Dai, Bulbul, ELIASZUNIGA2014849, Sayevand, Novin,  Zhang, Stoyanov, ELNAGGAR20161581, Weli, Hosen, CHOWDHURY20173962, Karahan}. Additionally, literatures describe the application of the cubic duffing equation, using Jacobi elliptic functions\cite{KOVACIC20161, ELIASZUNIGA20132574, Belendez}. As in the leaf function represented by Eq. (\ref{1.1.2}), the term $t^2$ is added to the root of inverse Jacobian elliptic function sn in Eq. (\ref{1.1.33}). Variables $r$ and $k$ control Jacobian function sn. The scope of applying the duffing equation to the Jacobian elliptic function is wider compared with the leaf function that has only one parameter $l$. Over time, the nonlinear duffing equation has witnessed improvements and further numerical analysis\cite{Nwamba, Serge, REMMI20182085, Koudahoun}.

\begin{equation}
\mathrm{Cubic-Quintic-Septic \: Duffing \: Equation:} \frac{\mathrm{d}^2r}{\mathrm{d}l^2}+ \alpha r+ \beta r^3+ \mu r^5+ \delta r^7=0 \label{1.1.36}
\end{equation}

A high-order exact solution using the Jacobi elliptic function has not been found yet. Furthermore, the high-order addition theorem necessary to derive the exact solution is not defined in the Jacobi elliptic function. To find the exact solution, it is key to determine if the superposition principle can be applied to linear equations. In the equations of various functions, mathematical operations divide one term into two or integrate two terms into one using the addition theorem. In this paper, we derive the addition theorem to further derive an exact solution for a high-order duffing equation, followed by applying the superposition principle.

%%%%%%%%% Section 2 %%%%%%%%%%%%%%%%%%%%%%%

\section{Addition Formulas}
\label{Addition Formulas}
\subsection{Addition Formulas of Leaf Functions}
\label{Addition Formulas of Leaf Functions}
Let there be two variables, $l_1$ and $l_2$. The addition formulas of the function $\mathrm{sleaf}_2(l)$ can be stated as follows:

\begin{equation}
\mathrm{sleaf}_2(l_1+l_2)=\frac{
\mathrm{sleaf}_2(l_1)\frac{\mathrm{\partial} \mathrm{sleaf}_2(l_2)}{\mathrm{\partial} l_2}
+\mathrm{sleaf}_2(l_2) \frac{\mathrm{\partial} \mathrm{sleaf}_2(l_1) }{\mathrm{\partial} l_1}
}
{1+(\mathrm{sleaf}_2(l_1))^2(\mathrm{sleaf}_2(l_2))^2}  \label{2.1.3}
\end{equation}

Depending on the domain of the variable $l$ of the leaf function, the signs of both $\mathrm{\partial} \mathrm{sleaf}_2(l_2)/\mathrm{\partial} l_2$ and $\mathrm{\partial} \mathrm{sleaf}_2(l_1)/\mathrm{\partial} l_1$ change. Eq. (\ref{2.1.3}) can be summarized  according a number of cases based on the domains of variables $l_1$ and $l_2$(See Figure 3). The parameters $m$ and $k$ represent integers.

(i) In the case where $\frac{\pi_2}{2}(4m-1) \leqq l_1 \leqq \frac{\pi_2}{2}(4m+1)$, $\frac{\pi_2}{2}(4k-1) \leqq l_2 \leqq \frac{\pi_2}{2}(4k+1)$(see Appendix E for the constant $\pi_2$), Eq. (\ref{2.1.3}) is transformed into:

\begin{equation}
\mathrm{sleaf}_2(l_1+l_2)=\frac{
\mathrm{sleaf}_2(l_1)\sqrt{1-(\mathrm{sleaf}_2(l_2))^4}
+\mathrm{sleaf}_2(l_2) \sqrt{1-(\mathrm{sleaf}_2(l_1))^4}
}
{1+(\mathrm{sleaf}_2(l_1))^2(\mathrm{sleaf}_2(l_2))^2}  \label{2.1.4}
\end{equation}

(ii) In the case where $\frac{\pi_2}{2}(4m-1) \leqq l_1 \leqq \frac{\pi_2}{2}(4m+1)$, $\frac{\pi_2}{2}(4k+1) \leqq l_2 \leqq \frac{\pi_2}{2}(4k+3)$, Eq. (\ref{2.1.3}) is transformed into:

\begin{equation}
\mathrm{sleaf}_2(l_1+l_2)=\frac{
-\mathrm{sleaf}_2(l_1)\sqrt{1-(\mathrm{sleaf}_2(l_2))^4}
+\mathrm{sleaf}_2(l_2) \sqrt{1-(\mathrm{sleaf}_2(l_1))^4}
}
{1+(\mathrm{sleaf}_2(l_1))^2(\mathrm{sleaf}_2(l_2))^2}   \label{2.1.5}
\end{equation}

(iii) In the case where $\frac{\pi_2}{2}(4m+1) \leqq l_1 \leqq \frac{\pi_2}{2}(4m+3)$, $\frac{\pi_2}{2}(4k-1) \leqq l_2 \leqq \frac{\pi_2}{2}(4k+1)$, Eq. (\ref{2.1.3}) is tansformed  into:

\begin{equation}
\mathrm{sleaf}_2(l_1+l_2)=\frac{
\mathrm{sleaf}_2(l_1)\sqrt{1-(\mathrm{sleaf}_2(l_2))^4}
-\mathrm{sleaf}_2(l_2) \sqrt{1-(\mathrm{sleaf}_2(l_1))^4}
}
{1+(\mathrm{sleaf}_2(l_1))^2(\mathrm{sleaf}_2(l_2))^2}   \label{2.1.6}
\end{equation}

(iiii) In the case where $\frac{\pi_2}{2}(4m+1) \leqq l_1 \leqq \frac{\pi_2}{2}(4m+3)$, $\frac{\pi_2}{2}(4k+1) \leqq l_2 \leqq \frac{\pi_2}{2}(4k+3)$, Eq. (\ref{2.1.3}) is tansformed  into:

\begin{equation}
\mathrm{sleaf}_2(l_1+l_2)=\frac{
-\mathrm{sleaf}_2(l_1)\sqrt{1-(\mathrm{sleaf}_2(l_2))^4}
-\mathrm{sleaf}_2(l_2) \sqrt{1-(\mathrm{sleaf}_2(l_1))^4}
}
{1+(\mathrm{sleaf}_2(l_1))^2(\mathrm{sleaf}_2(l_2))^2}   \label{2.1.7}
\end{equation}

Next, the addition formula of $\mathrm{cleaf}_2(l)$ can be stated as follows:

\begin{equation}
\mathrm{cleaf}_2(l_1+l_2)=\frac{
\mathrm{cleaf}_2(l_1)\frac{\mathrm{\partial sleaf}_2(l_2)}{\mathrm{\partial} l_2}
+\mathrm{sleaf}_2(l_2) \frac{\mathrm{\partial cleaf}_2(l_1) }{\mathrm{\partial} l_1}
}
{1+(\mathrm{cleaf}_2(l_1))^2(\mathrm{sleaf}_2(l_2))^2}   \label{2.1.8}
\end{equation}

Depending on the domain of the variable $l$ of the leaf function, the signs of both $\mathrm{\partial sleaf}_2(l_2)/\mathrm{\partial} l_2$ and $\mathrm{\partial cleaf}_2(l_1)/\mathrm{\partial} l_1$ change. Eq. (\ref{2.1.8}) can be summarized according a number of cases based on the domains of variables $l_1$ and $l_2$.

(i)  In the case where $2m \pi_2 \leqq l_1 \leqq (2m+1) \pi_2$, $\frac{\pi_2}{2}(4k-1) \leqq l_2 \leqq \frac{\pi_2}{2}(4k+1)$, Eq. (\ref{2.1.8}) is transformed into: 

\begin{equation}
\mathrm{cleaf}_2(l_1+l_2)=\frac{
\mathrm{cleaf}_2(l_1)\sqrt{1-(\mathrm{sleaf}_2(l_2))^4}
-\mathrm{sleaf}_2(l_2) \sqrt{1-(\mathrm{cleaf}_2(l_1))^4}
}
{1+(\mathrm{cleaf}_2(l_1))^2(\mathrm{sleaf}_2(l_2))^2}   \label{2.1.9}
\end{equation}

(ii) In the case where $(2m-1) \pi_2 \leqq l_1 \leqq 2m \pi_2$, $\frac{\pi_2}{2}(4k-1) \leqq l_2 \leqq \frac{\pi_2}{2}(4k+1)$, Eq. (\ref{2.1.8}) is transformed into: 

\begin{equation}
\mathrm{cleaf}_2(l_1+l_2)=\frac{
\mathrm{cleaf}_2(l_1)\sqrt{1-(\mathrm{sleaf}_2(l_2))^4}
+\mathrm{sleaf}_2(l_2) \sqrt{1-(\mathrm{cleaf}_2(l_1))^4}
}
{1+(\mathrm{cleaf}_2(l_1))^2(\mathrm{sleaf}_2(l_2))^2}   \label{2.1.10}
\end{equation}

(iii) In the case where $(2m-1) \pi_2 \leqq l_1 \leqq 2m \pi_2$, $\frac{\pi_2}{2}(4k+1) \leqq l_2 \leqq \frac{\pi_2}{2}(4k+3)$, Eq. (\ref{2.1.8}) is transformed into:

\begin{equation}
\mathrm{cleaf}_2(l_1+l_2)=\frac{
-\mathrm{cleaf}_2(l_1)\sqrt{1-(\mathrm{sleaf}_2(l_2))^4}
+\mathrm{sleaf}_2(l_2) \sqrt{1-(\mathrm{cleaf}_2(l_1))^4}
}
{1+(\mathrm{cleaf}_2(l_1))^2(\mathrm{sleaf}_2(l_2))^2}   \label{2.1.11}
\end{equation}

(iiii) In the case where $2m \pi_2 \leqq l_1 \leqq (2m+1) \pi_2$, $\frac{\pi_2}{2}(4k+1) \leqq l_2 \leqq \frac{\pi_2}{2}(4k+3)$, Eq. (\ref{2.1.8}) is transformed into:

\begin{equation}
\mathrm{cleaf}_2(l_1+l_2)=\frac{
-\mathrm{cleaf}_2(l_1)\sqrt{1-(\mathrm{sleaf}_2(l_2))^4}
-\mathrm{sleaf}_2(l_2) \sqrt{1-(\mathrm{cleaf}_2(l_1))^4}
}
{1+(\mathrm{cleaf}_2(l_1))^2(\mathrm{sleaf}_2(l_2))^2}  \label{2.1.12}
\end{equation}

Next, the addition formulas of $\mathrm{sleaf}_3(l)$ can be described as follows:

\begin{equation}
\begin{split}
(\mathrm{sleaf}_3(l_1+l_2))^2=\frac{ \left\{
\mathrm{sleaf}_3(l_1)\frac{\mathrm{\partial sleaf}_3(l_2)}{\mathrm{\partial} l_2}
+\mathrm{sleaf}_3(l_2) \frac{\mathrm{\partial sleaf}_3(l_1) }{\mathrm{\partial} l_1}
\right\}^2
}
{1+4(\mathrm{sleaf}_3(l_1))^4(\mathrm{sleaf}_3(l_2))^2+4(\mathrm{sleaf}_3(l_1))^2(\mathrm{sleaf}_3(l_2))^4} \\ 
+\frac{\left\{
(\mathrm{sleaf}_3(l_1))^3 \mathrm{sleaf}_3(l_2) - \mathrm{sleaf}_3(l_1)	(\mathrm{sleaf}_3(l_2))^3 
\right\}^2}
{1+4(\mathrm{sleaf}_3(l_1))^4(\mathrm{sleaf}_3(l_2))^2+4(\mathrm{sleaf}_3(l_1))^2(\mathrm{sleaf}_3(l_2))^4}
\label{2.1.13}
\end{split}
\end{equation}

The preceding equation can be summarized as follows according to a number of cases based on the domains of the variables $l_1$ and $l_2$.

(i) In the case where both $(4m-1) \frac{\pi_3}{2} \leqq l_1 
\leqq (4m+1) \frac{\pi_3}{2}$ and $(4k-1) \frac{\pi_3}{2} \leqq l_2 
\leqq (4k+1) \frac{\pi_3}{2}$ 
or both $(4m+1) \frac{\pi_3}{2} \leqq l_1 \leqq (4m+3) \frac{\pi_3}{2}$ and $(4k+1) \frac{\pi_3}{2} \leqq l_2 \leqq (4k+3) \frac{\pi_3}{2}$, Eq. (\ref{2.1.13}) is tansformed into: 

\begin{equation}
\begin{split}
(\mathrm{sleaf}_3(l_1+l_2))^2=\frac{ \left\{
\mathrm{sleaf}_3(l_1)\sqrt{1-(\mathrm{sleaf}_3(l_2))^6}
+\mathrm{sleaf}_3(l_2) \sqrt{1-(\mathrm{sleaf}_3(l_1))^6}
\right\}^2
}
{1+4(\mathrm{sleaf}_3(l_1))^4(\mathrm{sleaf}_3(l_2))^2+4(\mathrm{sleaf}_3(l_1))^2(\mathrm{sleaf}_3(l_2))^4} \\ 
+\frac{ \left\{
(\mathrm{sleaf}_3(l_1))^3 \mathrm{sleaf}_3(l_2) - \mathrm{sleaf}_3(l_1)	(\mathrm{sleaf}_3(l_2))^3 
\right\}^2}
{1+4(\mathrm{sleaf}_3(l_1))^4(\mathrm{sleaf}_3(l_2))^2+4(\mathrm{sleaf}_3(l_1))^2(\mathrm{sleaf}_3(l_2))^4}
\label{2.1.14}
\end{split}
\end{equation}

The symbol $\pi_3$ represents a constant (see Appendix E).

(ii) In the case where both $(4m+1) \frac{\pi_3}{2} \leqq l_1 
\leqq (4m+3) \frac{\pi_3}{2}$ and $(4k-1) \frac{\pi_3}{2} \leqq l_2 \leqq (4k+1) \frac{\pi_3}{2}$ 
or both $(4m-1) \frac{\pi_3}{2} \leqq l_1 \leqq (4m+1) \frac{\pi_3}{2}$ and $(4k+1) \frac{\pi_3}{2} \leqq l_2 \leqq (4k+3) \frac{\pi_3}{2}$, Eq. (\ref{2.1.13}) is tansformed into: 

\begin{equation}
\begin{split}
(\mathrm{sleaf}_3(l_1+l_2))^2=\frac{ \left\{
\mathrm{sleaf}_3(l_1)\sqrt{1-(\mathrm{sleaf}_3(l_2))^6}
-\mathrm{sleaf}_3(l_2) \sqrt{1-(\mathrm{sleaf}_3(l_1))^6}
\right\}^2
}
{1+4(\mathrm{sleaf}_3(l_1))^4(\mathrm{sleaf}_3(l_2))^2+4(\mathrm{sleaf}_3(l_1))^2(\mathrm{sleaf}_3(l_2))^4} \\ 
+\frac{ \left\{
(\mathrm{sleaf}_3(l_1))^3 \mathrm{sleaf}_3(l_2) - \mathrm{sleaf}_3(l_1)	(\mathrm{sleaf}_3(l_2))^3 
\right\}^2}
{1+4(\mathrm{sleaf}_3(l_1))^4(\mathrm{sleaf}_3(l_2))^2+4(\mathrm{sleaf}_3(l_1))^2(\mathrm{sleaf}_3(l_2))^4}
\label{2.1.15}
\end{split}
\end{equation}

Next, the addition formulas of $\mathrm{cleaf}_3(l)$ can be defined as follows:

\begin{equation}
\begin{split}
(\mathrm{cleaf}_3(l_1+l_2))^2=\frac{\left\{
\mathrm{cleaf}_3(l_1)\frac{\mathrm{\partial sleaf}_3(l_2)}{\mathrm{\partial} l_2}
+\mathrm{sleaf}_3(l_2) \frac{\mathrm{\partial cleaf}_3(l_1) }{\mathrm{\partial} l_1}
\right\}^2
}
{1+4(\mathrm{sleaf}_3(l_2))^4(\mathrm{cleaf}_3(l_1))^2+4(\mathrm{sleaf}_3(l_2))^2(\mathrm{cleaf}_3(l_1))^4} \\ 
+\frac{\left\{
(\mathrm{sleaf}_3(l_2))^3 \mathrm{cleaf}_3(l_1) - \mathrm{sleaf}_3(l_2)	(\mathrm{cleaf}_3(l_1))^3
\right\}^2}
{1+4(\mathrm{sleaf}_3(l_2))^4(\mathrm{cleaf}_3(l_1))^2+4(\mathrm{sleaf}_3(l_2))^2(\mathrm{cleaf}_3(l_1))^4}
\label{2.1.16}
\end{split}
\end{equation}

The preceding equation can be summarized as follows according to a number of cases based on the domains of the variables $l_1$ and $l_2$.

(i) In the case where both $2k \pi_3 \leqq l_1\leqq (2k+1) \pi_3$ and $(4m-1) \frac{\pi_3}{2} \leqq l_2 
\leqq (4m+1) \frac{\pi_3}{2}$ 
or both $(2k+1) \pi_3 \leqq l_1 \leqq (2k+2) \pi_3 $ and $(4m+1) \frac{\pi_3}{2} \leqq l_2 \leqq (4m+3) \frac{\pi_3}{2}$, Eq. (\ref{2.1.16}) is transformed into: 

\begin{equation}
\begin{split}
(\mathrm{cleaf}_3(l_1+l_2))^2=\frac{\left\{
\mathrm{cleaf}_3(l_1) \sqrt{1- (\mathrm{sleaf}_3(l_2))^6 }
-\mathrm{sleaf}_3(l_2) \sqrt{1- (\mathrm{cleaf}_3(l_1))^6 }
\right\}^2
}
{1+4(\mathrm{sleaf}_3(l_2))^4(\mathrm{cleaf}_3(l_1))^2+4(sleaf_3(l_2))^2(\mathrm{cleaf}_3(l_1))^4} \\ 
+\frac{\left\{
(\mathrm{sleaf}_3(l_2))^3 \mathrm{cleaf}_3(l_1) - \mathrm{sleaf}_3(l_2)	(\mathrm{cleaf}_3(l_1))^3 
\right\}^2}
{1+4(\mathrm{sleaf}_3(l_2))^4(\mathrm{cleaf}_3(l_1))^2+4(\mathrm{sleaf}_3(l_2))^2(\mathrm{cleaf}_3(l_1))^4
}
\label{2.1.17}
\end{split}
\end{equation}

(ii) In the case where both $(2k+1) \pi_3 \leqq l_1\leqq (2k+2) \pi_3$ and $(4m-1) \frac{\pi_3}{2} \leqq l_2 
\leqq (4m+1) \frac{\pi_3}{2}$ 
or both $2k \pi_3 \leqq l_1 \leqq (2k+1) \pi_3 $ and $(4m+1) \frac{\pi_3}{2} \leqq l_2 \leqq (4m+3) \frac{\pi_3}{2}$, Eq. (\ref{2.1.16}) is transformed into: 

\begin{equation}
\begin{split}
(\mathrm{cleaf}_3(l_1+l_2))^2=\frac{\left\{
\mathrm{cleaf}_3(l_1) \sqrt{1- (\mathrm{sleaf}_3(l_2))^6 }
+\mathrm{sleaf}_3(l_2) \sqrt{1- (\mathrm{cleaf}_3(l_1))^6 }
\right\}^2}
{1+4(\mathrm{sleaf}_3(l_2))^4(\mathrm{cleaf}_3(l_1))^2+4(\mathrm{sleaf}_3(l_2))^2(\mathrm{cleaf}_3(l_1))^4} \\ 
+\frac{\left\{
(\mathrm{sleaf}_3(l_2))^3 \mathrm{cleaf}_3(l_1) - \mathrm{sleaf}_3(l_2)	(\mathrm{cleaf}_3(l_1))^3
\right\}^2}
{1+4(\mathrm{sleaf}_3(l_2))^4(\mathrm{cleaf}_3(l_1))^2+4(\mathrm{sleaf}_3(l_2))^2((\mathrm{cleaf}_3(l_1))^4}
\label{2.1.18}
\end{split}
\end{equation}

\begin{figure*}[tb]
\begin{center} 
\includegraphics[width=0.75 \textwidth]{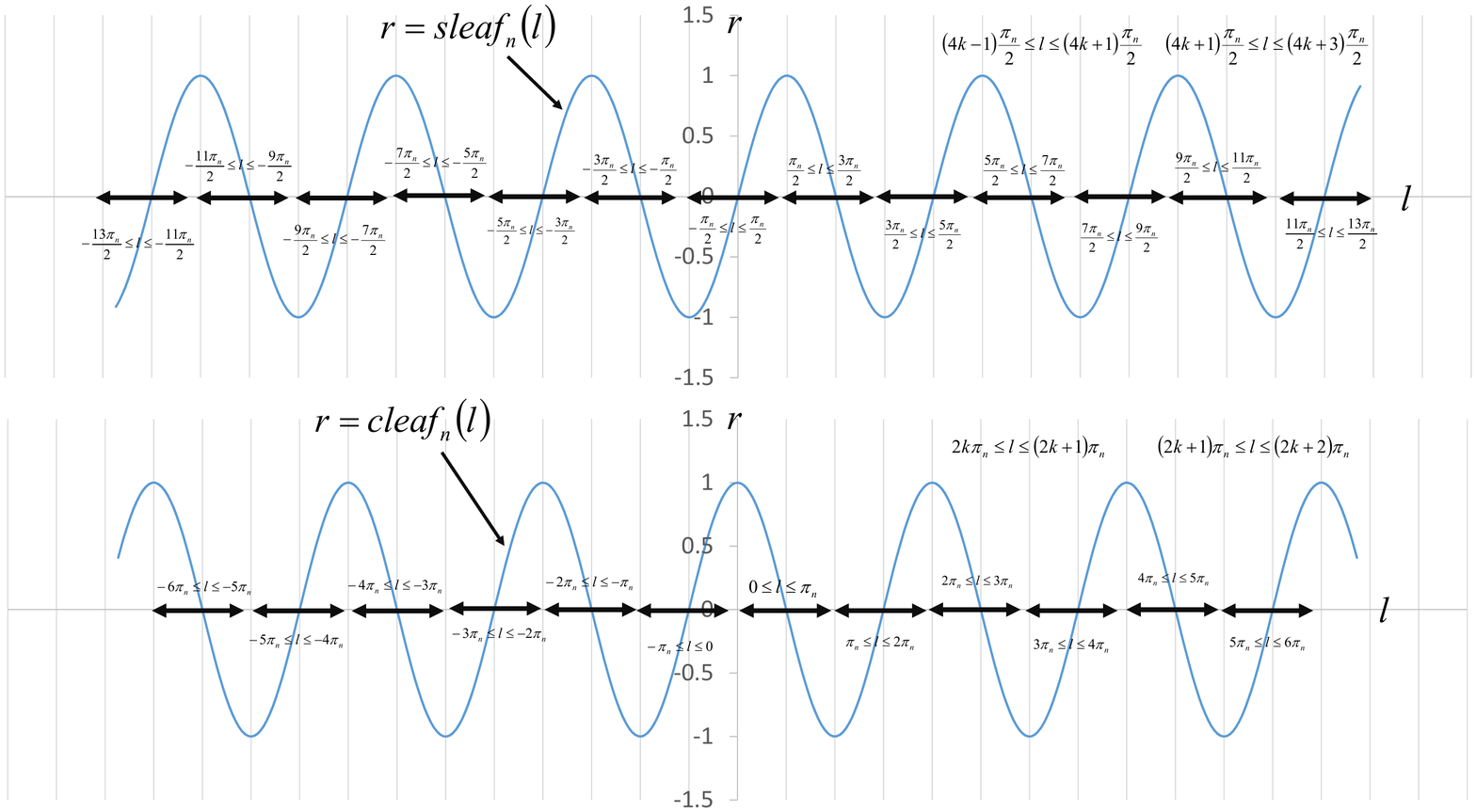}
\caption{ Curves of the functions $\mathrm{sleaf}_n(l)$ and $\mathrm{cleaf}_n(l)$ }
\label{fig2.1.1}
\end{center}
\end{figure*}

\subsection{Addition Formulas of Hyperbolic Leaf Function}
\label{Addition Formulas of Hyperbolic Leaf Function}

Let there be two variables, $l_1$ and $l_2$. Considering the imaginary number $i$, the relation between $\mathrm{sleaf}_2(l)$ and $\mathrm{sleafh}_2(l)$, and the relation between $\mathrm{cleaf}_2(l)$ and $\mathrm{cleafh}_2(l)$ can be obtained as follows(see Appendix D in detail):

\begin{equation}
\mathrm{sleaf}_2(i \cdot l)=i \cdot \mathrm{sleaf}_2(l)  \label{2.2.7}
\end{equation}

\begin{equation}
\mathrm{sleafh}_2( i \cdot l)=i \cdot \mathrm{sleafh}_2(l)  \label{2.2.8}
\end{equation}

\begin{equation}
\mathrm{cleaf}_2(i \cdot l)=\mathrm{cleafh}_2(l)  \label{2.2.9}
\end{equation}

As shown in Eq. (\ref{2.2.7}) and Eq. (\ref{2.2.8}), in the case where $n = 2$, 
the functions $\mathrm{sleaf}_2(i \cdot l)$ and $\mathrm{sleafh}_2(i \cdot l)$ are equal to the functions  $i \cdot \mathrm{sleaf}_2(l)$ and $i \cdot \mathrm{sleafh}_2(l)$, respectively.
Therefore, we cannot derive the addition formulas of $\mathrm{sleafh}_2(l)$ by replacing $i \cdot l$ with $l$ in Eqs. (\ref{2.1.4}) $\sim $ (\ref{2.1.7}). Using the relation between the function $\mathrm{sleaf}_2(l)$ and the function $\mathrm{sleafh}_2(l)$(see Appendix B),  the addition formulas of $\mathrm{sleafh}_2(l)$ can be obtained.
By substituting Eq. (\ref{B4}) into Eqs. (\ref{2.1.4}) $\sim $ (\ref{2.1.7}), the following equation is obtained:

\begin{equation}
\begin{split}
&\mathrm{sleafh}_2(l_1+l_2)= \\
&\frac{
\mathrm{sleafh}_2(l_1)\sqrt{1+(\mathrm{sleafh}_2(l_2))^4}
+\mathrm{sleafh}_2(l_2) \sqrt{1+(\mathrm{sleafh}_2(l_1))^4}
}
{1-(\mathrm{sleafh}_2(l_1))^2(\mathrm{sleafh}_2(l_2))^2}  \label{2.2.10}
\end{split}
\end{equation}

In the work\cite{Kaz_ch}, the addition formulas of $\mathrm{cleafh}_2(l)$ are obtained using Eq. (\ref{B3}):

\begin{equation}
\begin{split}
&\mathrm{cleafh}_2(l_1+l_2)= \\
&\frac{
2 \mathrm{cleafh}_2(l_1) \mathrm{cleafh}_2(l_2) 
+\frac{\mathrm{\partial cleaf}_2(l_1)}{ \mathrm{\partial} l_1}
\frac{ \mathrm{\partial cleaf}_2(l_2)}{ \mathrm{\partial} l_2}
}
{1+(\mathrm{cleafh}_2(l_1))^2+(\mathrm{cleafh}_2(l_2))^2-(\mathrm{cleafh}_2(l_1))^2 (\mathrm{cleafh}_2(l_2))^2}  \label{2.2.11}
\end{split}
\end{equation}

The preceding equation can be summarized as follows according to a number of cases based on the domains of variables $l_1$ and $l_2$.

(i) In the case where both $0 \leqq l_1 \leqq \eta_2$ and $0 \leqq l_2 \leqq \eta_2$ or both $-\eta_2 \leqq l_1 \leqq 0$ and $-\eta_2 \leqq l_2 \leqq 0$ (see Appendix G for the constant $\eta_2$.), Eq. (\ref{2.2.11}) is transformed into:

\begin{equation}
\begin{split}
&\mathrm{cleafh}_2(l_1+l_2)= \\
&\frac{
2\mathrm{cleafh}_2(l_1) \mathrm{cleafh}_2(l_2) 
+\sqrt{ (\mathrm{cleafh}_2(l_1))^4-1 } \sqrt{ (\mathrm{cleafh}_2(l_2))^4-1 }
}
{1+(\mathrm{cleafh}_2(l_1))^2+(\mathrm{cleafh}_2(l_2))^2-(\mathrm{cleafh}_2(l_1))^2 (\mathrm{cleafh}_2(l_2))^2}  \label{2.2.12}
\end{split}
\end{equation}

(ii)  In the case where both $0 \leqq l_1 \leqq \eta_2$ and $-\eta_2 \leqq l_2 \leqq 0$ or both $-\eta_2 \leqq l_1 \leqq 0$ and $0 \leqq l_2 \leqq  \eta_2$, Eq. (\ref{2.2.11}) is transformed into:

\begin{equation}
\begin{split}
&\mathrm{cleafh}_2(l_1+l_2)= \\
&\frac{
2\mathrm{cleafh}_2(l_1) \mathrm{cleafh}_2(l_2) 
-\sqrt{ (\mathrm{cleafh}_2(l_1))^4-1 } \sqrt{ (\mathrm{cleafh}_2(l_2))^4-1 }
}
{1+(\mathrm{cleafh}_2(l_1))^2+(\mathrm{cleafh}_2(l_2))^2-(\mathrm{cleafh}_2(l_1))^2 (\mathrm{cleafh}_2(l_2))^2}  \label{2.2.13}
\end{split}
\end{equation}

Next, let us consider the case of $n=3$. The relation between $\mathrm{sleaf}_3(l)$ and $\mathrm{sleafh}_3(l)$ and the relation between $\mathrm{cleaf}_3(l)$ and $\mathrm{cleafh}_3(l)$ are as follows (see Appendix D):

\begin{equation}
\mathrm{sleaf}_3(l)=-i \cdot \mathrm{sleafh}_3(i \cdot l)   \label{2.2.14}
\end{equation}

\begin{equation}
\mathrm{cleaf}_3(l)=\mathrm{cleafh}_3(i \cdot l)   \label{2.2.15}
\end{equation}

In Eq. (\ref{2.1.13}) $\sim$ Eq. (\ref{2.1.18}), the variables $l_1$ and $l_2$ are replaced with the expressions $i \cdot l_1$ and $i \cdot l_2$, respectively. The addition formulas of $\mathrm{sleafh}_3(l)$ are defined as follows:

\begin{equation}
\begin{split}
& (\mathrm{sleafh}_3(l_1+l_2))^2= \\ 
& \frac{\left\{
\mathrm{sleafh}_3(l_1)\sqrt{1+(\mathrm{sleafh}_3(l_2))^6}+\mathrm{sleafh}_3(l_2) \sqrt{1+(\mathrm{sleafh}_3(l_1))^6}
\right\}^2}
{1-4(\mathrm{sleafh}_3(l_1))^4(\mathrm{sleafh}_3(l_2))^2-4(\mathrm{sleafh}_3(l_1))^2(\mathrm{sleafh}_3(l_2))^4}  \\
& -\frac{\left\{
(\mathrm{sleafh}_3(l_1))^3 \mathrm{sleafh}_3(l_2) - \mathrm{sleafh}_3(l_1)	(\mathrm{sleafh}_3(l_2))^3
\right\}^2}
{1-4(\mathrm{sleafh}_3(l_1))^4(\mathrm{sleafh}_3(l_2))^2-4(\mathrm{sleafh}_3(l_1))^2(\mathrm{sleafh}_3(l_2))^4}
 \label{2.2.16}
\end{split}
\end{equation}

The addition formulas of $\mathrm{cleafh}_3(l)$ are defined as follows:

\begin{equation}
\begin{split}
& (\mathrm{cleafh}_3(l_1+l_2))^2= \\ 
& \frac{\left\{
\mathrm{cleafh}_3(l_1) \frac{\mathrm{\partial sleafh}_3(l_2)}{\mathrm{\partial} l_2}
+\mathrm{sleafh}_3(l_2) \frac{\mathrm{\partial cleafh}_3(l_1)}{\mathrm{\partial} l_1}
\right\}^2}
{1+4(\mathrm{sleafh}_3(l_2))^4(\mathrm{cleafh}_3(l_1))^2-4(\mathrm{sleafh}_3(l_2))^2(\mathrm{cleafh}_3(l_1))^4}  \\
& -\frac{\left\{
(\mathrm{sleafh}_3(l_1))^3 \mathrm{cleafh}_3(l_2) + \mathrm{sleafh}_3(l_2)	(\mathrm{cleafh}_3(l_1))^3
\right\}^2}
{1+4(\mathrm{sleafh}_3(l_2))^4(\mathrm{cleafh}_3(l_1))^2-4(\mathrm{sleafh}_3(l_2))^2(\mathrm{cleafh}_3(l_1))^4}
 \label{2.2.17}
\end{split}
\end{equation}

The preceding equation can be summarized as follows according to a number of cases based on the domains of the variables $l_1$ and $l_2$.

(i)  In the case where both $-\eta_3 \leqq l_1 \leqq 0$  (see Appendix G for the constant $\eta_3$), Eq. (\ref{2.2.17}) is transformed into:

\begin{equation}
\begin{split}
& (\mathrm{cleafh}_3(l_1+l_2))^2= \\ 
& \frac{\left\{
\mathrm{cleafh}_3(l_1) \sqrt{1 + (\mathrm{sleafh}_3(l_2))^6 } 
-\mathrm{sleafh}_3(l_2) \sqrt{ (\mathrm{cleafh}_3(l_1))^6-1 }
\right\}^2}
{1+4(\mathrm{sleafh}_3(l_2))^4(\mathrm{cleafh}_3(l_1))^2-4(\mathrm{sleafh}_3(l_2))^2(\mathrm{cleafh}_3(l_1))^4}  \\
& -\frac{\left\{
(\mathrm{sleafh}_3(l_2))^3 \mathrm{cleafh}_3(l_1) + \mathrm{sleafh}_3(l_2)	(\mathrm{cleafh}_3(l_1))^3
\right\}^2}
{1+4(\mathrm{sleafh}_3(l_2))^4(\mathrm{cleafh}_3(l_1))^2-4(\mathrm{sleafh}_3(l_2))^2(\mathrm{cleafh}_3(l_1))^4}
 \label{2.2.18}
\end{split}
\end{equation}

(ii)  In the case where $0 \leqq l_1 \leqq \eta_3$, Eq. (\ref{2.2.17}) is transformed into:

\begin{equation}
\begin{split}
& (\mathrm{cleafh}_3(l_1+l_2))^2= \\ 
& \frac{\left\{
\mathrm{cleafh}_3(l_1) \sqrt{1 + (\mathrm{sleafh}_3(l_2))^6 } 
+\mathrm{sleafh}_3(l_2) \sqrt{ (\mathrm{cleafh}_3(l_1))^6-1 }
\right\}^2}
{1+4(\mathrm{sleafh}_3(l_2))^4(\mathrm{cleafh}_3(l_1))^2-4(\mathrm{sleafh}_3(l_2))^2(\mathrm{cleafh}_3(l_1))^4}  \\
& -\frac{\left\{
(\mathrm{sleafh}_3(l_2))^3 \mathrm{cleafh}_3(l_1) + \mathrm{sleafh}_3(l_2)	(\mathrm{cleafh}_3(l_1))^3
\right\}^2}
{1+4(\mathrm{sleafh}_3(l_2))^4(\mathrm{cleafh}_3(l_1))^2-4(\mathrm{sleafh}_3(l_2))^2(\mathrm{cleafh}_3(l_1))^4}
 \label{2.2.19}
\end{split}
\end{equation}

%%%%%%%%% Section 3 %%%%%%%%%%%%%%%%%%%%%%%

\section{Double Angle Formulas and Half Angle Formulas }
\label{Double Angle Formulas and Half Angle Formulass}
\subsection{Double Angle Formulas of Leaf Functions}
\label{Double Angle Formulas of Leaf Functions}
In the case where the basis $n=2$, the variables $l_1$ and $l_2$ in Eq. (\ref{2.1.3}) are replaced with the variable $l$, and the double-angle formula can be expressed as follows:

\begin{equation}
\mathrm{sleaf}_2(2l)=\frac{ 2\mathrm{sleaf}_2(l) \frac{\mathrm{\partial sleaf}_2(l)}{\mathrm{\partial} l} }{1+(\mathrm{sleaf}_2(l))^4}  \label{3.1.3}
\end{equation}

The preceding equation can be summarized as follows according to a number of cases based on the domain of the variable $l$.

(i) In the case where $ \frac{\pi_2}{2}(4m-1) \leqq l \leqq \frac{\pi_2}{2}(4m+1) $, Eq. (\ref{3.1.3}) is transformed into:

\begin{equation}
\mathrm{sleaf}_2(2l)=\frac{ 2\mathrm{sleaf}_2(l) \sqrt{1- (\mathrm{sleaf}_2(l))^4}}{1+(\mathrm{sleaf}_2(l))^4}  \label{3.1.4}
\end{equation}

(ii) In the case where $ \frac{\pi_2}{2}(4m+1) \leqq l \leqq \frac{\pi_2}{2}(4m+3) $, Eq. (\ref{3.1.3}) is transformed into:

\begin{equation}
\mathrm{sleaf}_2(2l)=-\frac{ 2\mathrm{sleaf}_2(l) \sqrt{1- (\mathrm{sleaf}_2(l))^4}}{1+(\mathrm{sleaf}_2(l))^4}  \label{3.1.5}
\end{equation}

The variables $l_1$ and $l_2$ in Eq. (\ref{2.1.9}) $\sim$ (\ref{2.1.12}) are replaced with the variable $l$. The double-angle formula can be defined as follows:

\begin{equation}
\mathrm{cleaf}_2(2l)=\frac{1-2 (\mathrm{cleaf}_2(l))^2-(\mathrm{cleaf}_2(l))^4 }
{-1-2 (\mathrm{cleaf}_2(l))^2+(\mathrm{cleaf}_2(l))^4 }  \label{3.1.6}
\end{equation}

In the case where the basis $n = 3$, the variables $l_1$ and $l_2$ of Eq. (\ref{2.1.13}) are replaced with the variable $l$, and the double-angle formula of the function $\mathrm{sleaf}_3(2l)$ can be expressed as follows:

\begin{equation}
\mathrm{sleaf}_3(2l)=\frac{2 \mathrm{sleaf}_3(l) \frac{\mathrm{\partial sleaf}_3(l)}{\partial l} }
{ \sqrt{1+8( \mathrm{sleaf}_3(l) )^6} }  \label{3.1.7}
\end{equation}

(i) In the case where $ \frac{\pi_3}{2}(4m-1) \leqq l \leqq \frac{\pi_3}{2}(4m+1) $ (see Appendix E for the constant $\pi_3$), Eq. (\ref{3.1.7}) is transformed into:

\begin{equation}
\mathrm{sleaf}_3(2l)=\frac{2 \mathrm{sleaf}_3(l) \sqrt{1 - (\mathrm{sleaf}_3(l))^6 } }
{ \sqrt{1+8( \mathrm{sleaf}_3(l))^6 } }  \label{3.1.8}
\end{equation}

(ii) In the case where $ \frac{\pi_3}{2}(4m+1) \leqq l \leqq \frac{\pi_3}{2}(4m+3) $, Eq. (\ref{3.1.7}) is transformed into:

\begin{equation}
\mathrm{sleaf}_3(2l)=-\frac{2 \mathrm{sleaf}_3(l) \sqrt{1 - (\mathrm{sleaf}_3(l))^6 } }
{ \sqrt{1+8( \mathrm{sleaf}_3(l))^6 } }  \label{3.1.9}
\end{equation}

In the case where the basis $n = 3$, the variable $l_1$ and the variable $l_2$ of Eq. (\ref{2.1.17}) $\sim$ Eq. (\ref{2.1.18}) are replaced with the variable $l$. The double-angle formula of the function $\mathrm{cleaf}_3(2l)$ is then expressed as follows:

\begin{equation}
\mathrm{cleaf}_3(2l)=\frac{2 (\mathrm{cleaf}_3(l))^2+2(\mathrm{cleaf}_3(l))^4-1 }
{ \sqrt{1+8( \mathrm{cleaf}_3(l))^2+8( \mathrm{cleaf}_3(l))^6-8( \mathrm{cleaf}_3(l))^8  } } \label{3.1.10}
\end{equation}

\subsection{Half Angle Formulas of Leaf Functions }
\label{Half Angle Formulas of Leaf Functions }
In the case where the basis $n=2$, the variables $l_1$ and $l_2$ in Eqs. (\ref{2.1.4}) $\sim$ (\ref{2.1.7}) are replaced with the expression $l/2$, and the half-angle formula is defined as follows:

(i) In the case where $ \frac{\pi_2}{2}(4m+1) \leqq l \leqq \frac{\pi_2}{2}(4m+3) $ (see Appendix E for the constant $\pi_2$), the half-angle formula is expressed as follows: 

\begin{equation}
\begin{split}
&\Bigl(\mathrm{sleaf}_2 \Bigl(\frac{l}{2} \Bigr) \Bigr)^2 = \frac{-1 - \sqrt{1 - (\mathrm{sleaf}_2(l))^2 } } { ( \mathrm{sleaf}_2(l) )^2 } \\
&+\frac{ \sqrt{1 + (\mathrm{sleaf}_2(l))^2 } }{ ( \mathrm{sleaf}_2(l))^2 }  
\sqrt{2 - (\mathrm{sleaf}_2(l))^2+2\sqrt{1 - (\mathrm{sleaf}_2(l))^2 } } 
\label{3.2.3}
\end{split}
\end{equation}

(ii) In the case where $ \frac{\pi_2}{2}(4m-1) \leqq l \leqq \frac{\pi_2}{2}(4m+1) $, the half-angle formula is defined as follows: 

\begin{equation}
\begin{split}
&\Bigl(\mathrm{sleaf}_2 \Bigl(\frac{l}{2} \Bigr) \Bigr)^2 = \frac{-1 + \sqrt{1 - (\mathrm{sleaf}_2(l))^2 } } { ( \mathrm{sleaf}_2(l) )^2 } \\
&+\frac{ \sqrt{1 + (\mathrm{sleaf}_2(l))^2 } }{ ( \mathrm{sleaf}_2(l))^2 }  
\sqrt{2 - (\mathrm{sleaf}_2(l))^2-2\sqrt{1 - (\mathrm{sleaf}_2(l))^2 } } 
\label{3.2.4}
\end{split}
\end{equation}

In the case where the basis $n=2$, the variables $l_1$ and $l_2$ in Eqs. (\ref{2.1.9})-(\ref{2.1.12}) are replaced with the expression $l/2$ and the half-angle formula is expressed as follows:

\begin{equation}
\Bigl(\mathrm{cleaf}_2 \Bigl(\frac{l}{2} \Bigr) \Bigr)^2 = \frac{-1 + \mathrm{cleaf}_2(l) +\sqrt{2} \sqrt{1+(\mathrm{cleaf}_2(l))^2 } }
 {1+  \mathrm{cleaf}_2(l)  } 
\label{3.2.5}
\end{equation}

In the case where the basis $n=3$, the variables $l_1$ and $l_2$ in Eqs. ((\ref{2.1.13})-((\ref{2.1.15}) are replaced with the expression $l/2$ and the half-angle formula of the function $\mathrm{sleaf}_3(l)$ is defined as follows:

(i) In the case where $ \frac{\pi_3}{2}(4m-1) \leqq l \leqq \frac{\pi_3}{2}(4m+1) $ (see Appendix E for the constant $\pi_3$), the half-angle formula is defined as follows:

\begin{equation}
\begin{split}
&\Bigl( \mathrm{sleaf}_3 \Bigl(\frac{l}{2} \Bigr) \Bigr)^2=-\frac{1}{2} (\mathrm{sleaf}_3(l))^2
+\frac{1}{2} \sqrt{1+(\mathrm{sleaf}_3(l))^2+(\mathrm{sleaf}_3(l))^4 } \\
&-\frac{1}{2} \sqrt{-1-(\mathrm{sleaf}_3(l))^2+2(\mathrm{sleaf}_3(l))^4
+\frac{2-2 (\mathrm{sleaf}_3(l))^6}{\sqrt{1+(\mathrm{sleaf}_3(l))^2+(\mathrm{sleaf}_3(l))^4}}
}
\label{3.2.6}
\end{split}
\end{equation}

(ii) In the case where $ \frac{\pi_3}{2}(4m+1) \leqq l \leqq \frac{\pi_3}{2}(4m+3) $, the half-angle formula is expressed as follows:

\begin{equation}
\begin{split}
&\Bigl( \mathrm{sleaf}_3 \Bigl(\frac{l}{2} \Bigr) \Bigr)^2=-\frac{1}{2} (\mathrm{sleaf}_3(l))^2
+\frac{1}{2} \sqrt{1+(\mathrm{sleaf}_3(l))^2+(\mathrm{sleaf}_3(l))^4 } \\
&+\frac{1}{2} \sqrt{-1-(\mathrm{sleaf}_3(l))^2+2(\mathrm{sleaf}_3(l))^4
+\frac{2-2 (\mathrm{sleaf}_3(l))^6}{\sqrt{1+(\mathrm{sleaf}_3(l))^2+(\mathrm{sleaf}_3(l))^4}}
}
\label{3.2.7}
\end{split}
\end{equation}

In the case where the basis $n=3$, the variables $l_1$ and $l_2$ in Eqs. (\ref{2.1.17})-(\ref{2.1.18}) are replaced with the expression $l/2$ and the half-angle formula of the function $\mathrm{cleaf}_3(\frac{l}{2})$ is defined as follows:

\begin{equation}
\begin{split}
&\Bigl(\mathrm{cleaf}_3 \Bigl(\frac{l}{2} \Bigr)\Bigr)^2=\frac{(\mathrm{cleaf}_3(l))^2-1}{4(\mathrm{cleaf}_3(l))^2+2}
+\frac{\sqrt{3}\sqrt{1+(\mathrm{cleaf}_3(l))^2+(\mathrm{cleaf}_3(l))^4} }
{2\sqrt{1+4(\mathrm{cleaf}_3(l))^2+4(\mathrm{cleaf}_3(l))^4} } \\
& + \frac{\sqrt{3} \mathrm{cleaf}_3(l) 
\sqrt{
-3-6(\mathrm{cleaf}_3(l))^2+2\sqrt{3} \{1+2(\mathrm{cleaf}_3(l))^2 \}
\sqrt{1+(\mathrm{cleaf}_3(l))^2+(\mathrm{cleaf}_3(l))^4}
}
}{ 2 \{1+2(\mathrm{cleaf}_3(l))^2 \}^{\frac{3}{2}} }
\label{3.2.8}
\end{split}
\end{equation}

\subsection{Double Angle Formulas of Hyperbolic Leaf Functions }
\label{Double-angle formulas of the hyperbolic leaf functions}

In the case where the basis $n=2$, the variables $l_1$ and $l_2$ in Eq. (\ref{2.2.10}) are replaced with the variable $l$, and the double-angle formula is defined as follows:

\begin{equation}
\mathrm{sleafh}_2(2l)=\frac{ 2\mathrm{sleafh}_2(l) \sqrt{ 1+(\mathrm{sleafh}_2(l))^4 }}
{1-(\mathrm{sleafh}_2(l))^4}  \label{3.3.3}
\end{equation}

The variables $l_1$ and $l_2$ in Eq. (\ref{2.2.12}) and Eq. (\ref{2.2.13})  are replaced with the variable $l$. The double-angle formula is then defined as follows:

\begin{equation}
\mathrm{cleafh}_2(2l)=\frac{ (\mathrm{cleafh}_2(l))^4+2(\mathrm{cleafh}_2(l))^2-1 }
{ -(\mathrm{cleafh}_2(l))^4+2(\mathrm{cleafh}_2(l))^2 + 1 }  \label{3.3.4}
\end{equation}

In the case where the basis $n = 3$, the variables $l_1$ and $l_2$ of Eq. (\ref{2.2.16}) are replaced with the variable $l$, and the double-angle formula of the function $\mathrm{sleafh}_3(2l)$ is defined as follows:

\begin{equation}
\mathrm{sleafh}_3(2l)=\frac{2 \mathrm{sleafh}_3(l) \sqrt{1 + (\mathrm{sleafh}_3(l))^6} }{  \sqrt{1-8( \mathrm{sleafh}_3(l) )^6} }  \label{3.3.5}
\end{equation}

In the case where the basis $n = 3$, the variables $l_1$ and $l_2$ of Eq. (\ref{2.2.18}) and Eq. (\ref{2.2.19}) are replaced with the variable $l$, and the double-angle formula of the function $\mathrm{cleafh}_3(2l)$ is defined as follows:

\begin{equation}
\mathrm{cleafh}_3(2l)=\frac{2 (\mathrm{cleafh}_3(l))^2+2(\mathrm{cleafh}_3(l))^4-1 }
{ \sqrt{1+8( \mathrm{cleafh}_3(l))^2+8( \mathrm{cleafh}_3(l))^6-8( \mathrm{cleafh}_3(l))^8  } }  \label{3.3.6}
\end{equation}

\subsection{Half Angle Formulas of Hyperbolic Leaf Functions}
\label{Half-angle formulas of leaf function of hyperbolic leaf function}
In the case where the basis $n=2$, the variables $l_1$ and $l_2$ in Eq. (\ref{2.2.10}) are replaced with the expression $l/2$, and the half-angle formula is defined as follows:

(i) In the case where $ |l| \leqq \zeta_2 $ (see Appendix F for the constant $\zeta_2$ and Appendix H for the periodicity $n=2$), the half-angle formula is expressed as follows:

\begin{equation}
\Bigl(\mathrm{sleafh}_2 \Bigl(\frac{l}{2} \Bigr)\Bigr)^2 = \frac{1 + \sqrt{1 + (\mathrm{sleafh}_2(l))^4 } }
 { ( \mathrm{sleafh}_2(l) )^2 } 
-\frac{ \sqrt{2} }{ \sqrt{ -1 + \sqrt{1+(\mathrm{sleafh}_2(l))^4}}}   
\label{3.4.3}
\end{equation}

(ii) In the case where $ \zeta_2 \leqq |l| $, the half-angle formula is defined as follows:

\begin{equation}
\Bigl(\mathrm{sleafh}_2 \Bigl(\frac{l}{2} \Bigr)\Bigr)^2 = \frac{1 + \sqrt{1 + (\mathrm{sleafh}_2(l))^4 } }
 { ( \mathrm{sleafh}_2(l) )^2 } 
+\frac{ \sqrt{2} }{ \sqrt{ -1 + \sqrt{1+(\mathrm{sleafh}_2(l))^4}}}   
\label{3.4.4}
\end{equation}

In the case where the basis $n=2$, the variables $l_1$ and $l_2$ in  Eq. (\ref{2.2.12}) and Eq. (\ref{2.2.13}) are replaced with the expression $l/2$, and the half-angle formula can be expressed as follows (see Appendix G for the constant $\eta_2$ and Appendix H for the periodicity $n=2$):

(i) In the case where $  |l| \leqq \eta_2  $,  the half-angle formula is defined as follows:

\begin{equation}
\Bigl(\mathrm{cleafh}_2 \Bigl(\frac{l}{2} \Bigr) \Bigr)^2 = \frac{-1 + \mathrm{cleafh}_2(l) +\sqrt{2} \sqrt{1+(\mathrm{cleafh}_2(l))^2 } } {1+  \mathrm{cleafh}_2(l)  } 
\label{3.4.5}
\end{equation}

(ii) In the case where $  \eta_2  \leqq |l|  $,  the half-angle formula is defined as follows:

\begin{equation}
\Bigl(\mathrm{cleafh}_2 \Bigl(\frac{l}{2} \Bigr) \Bigr)^2 = \frac{-1 + \mathrm{cleafh}_2(l) -\sqrt{2} \sqrt{1+(\mathrm{cleafh}_2(l))^2 } } {1+  \mathrm{cleafh}_2(l)  } 
\label{3.4.6}
\end{equation}

In the case where the basis $n=3$, the variables $l_1$ and $l_2$ in  Eq. (\ref{2.2.16}) are replaced with the expression $l/2$, and the half-angle formula of the function $\mathrm{sleafh}_3(l)$ is defined as follows:

\begin{equation}
\begin{split}
&\Bigl(\mathrm{sleafh}_3 \Bigl(\frac{l}{2} \Bigr) \Bigr)^2=-\frac{1}{2} (\mathrm{sleafh}_3(l))^2
-\frac{1}{2} \sqrt{1-(\mathrm{sleafh}_3(l))^2+(\mathrm{sleafh}_3(l))^4 } \\
&+\frac{1}{2} \sqrt{-1+(\mathrm{sleafh}_3(l))^2+2(\mathrm{sleafh}_3(l))^4
+\frac{2+2 (\mathrm{sleafh}_3(l))^6}{\sqrt{1-(\mathrm{sleafh}_3(l))^2+(\mathrm{sleafh}_3(l))^4}}
}
\label{3.4.7}
\end{split}
\end{equation}

In the case where the basis $n=3$, the variables $l_1$ and $l_2$ in Eq. (\ref{2.2.18}) and  Eq. (\ref{2.2.19}) are replaced with the expression $l/2$ and the half-angle formula of the function $\mathrm{cleafh}_3(l)$ is defined as follows:

\begin{equation}
\begin{split}
&\Bigl(\mathrm{cleafh}_3 \Bigl(\frac{l}{2} \Bigr) \Bigr)^2=
\frac{-1+\mathrm{cleafh}_3(l))^2+\sqrt{3}\sqrt{1+\mathrm{cleafh}_3(l))^2+(\mathrm{cleafh}_3(l))^4}}{4(\mathrm{cleafh}_3(l))^2+2} \\
&+\frac{
\sqrt{3} \mathrm{cleafh}_3(l)
\sqrt{-3-6(\mathrm{cleafh}_3(l))^2+2\sqrt{3} \{ 1+2(\mathrm{cleafh}_3(l))^2 \}
\sqrt{1+(\mathrm{cleafh}_3(l))^2+(\mathrm{cleafh}_3(l))^4}
}
}
{
2 \{1+2 \mathrm{cleafh}_3(l))^2 \}^{3/2}
}
\label{3.4.8}
\end{split}
\end{equation}

\section{Numerical Analysis}
\label{Numerical Analysis}
\subsection{Numerical Analysis of Leaf Functions}
\label{Numerical Analysis of Leaf Functions}
\begin{figure*}[tb]
\begin{center} 
\includegraphics[width=0.75 \textwidth]{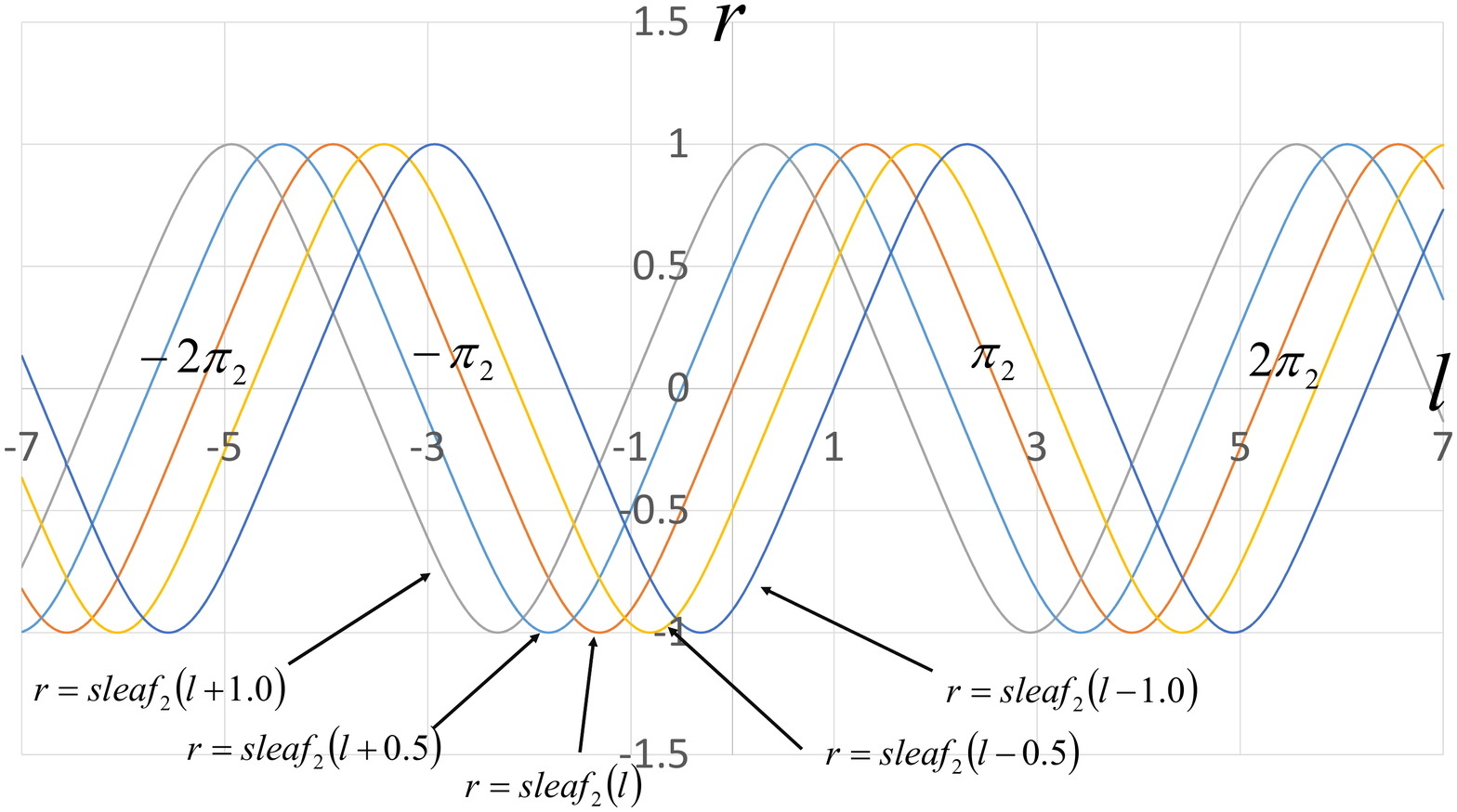}
\caption{ Translation of the curves of the function $\mathrm{sleaf}_2(l)$ obtained using the addition formulas with the basis $n = 2$ }
\label{fig4.1.4}        % Give a unique label
\end{center}
\end{figure*}

\begin{figure*}[tb]
\begin{center}
\includegraphics[width=0.75 \textwidth]{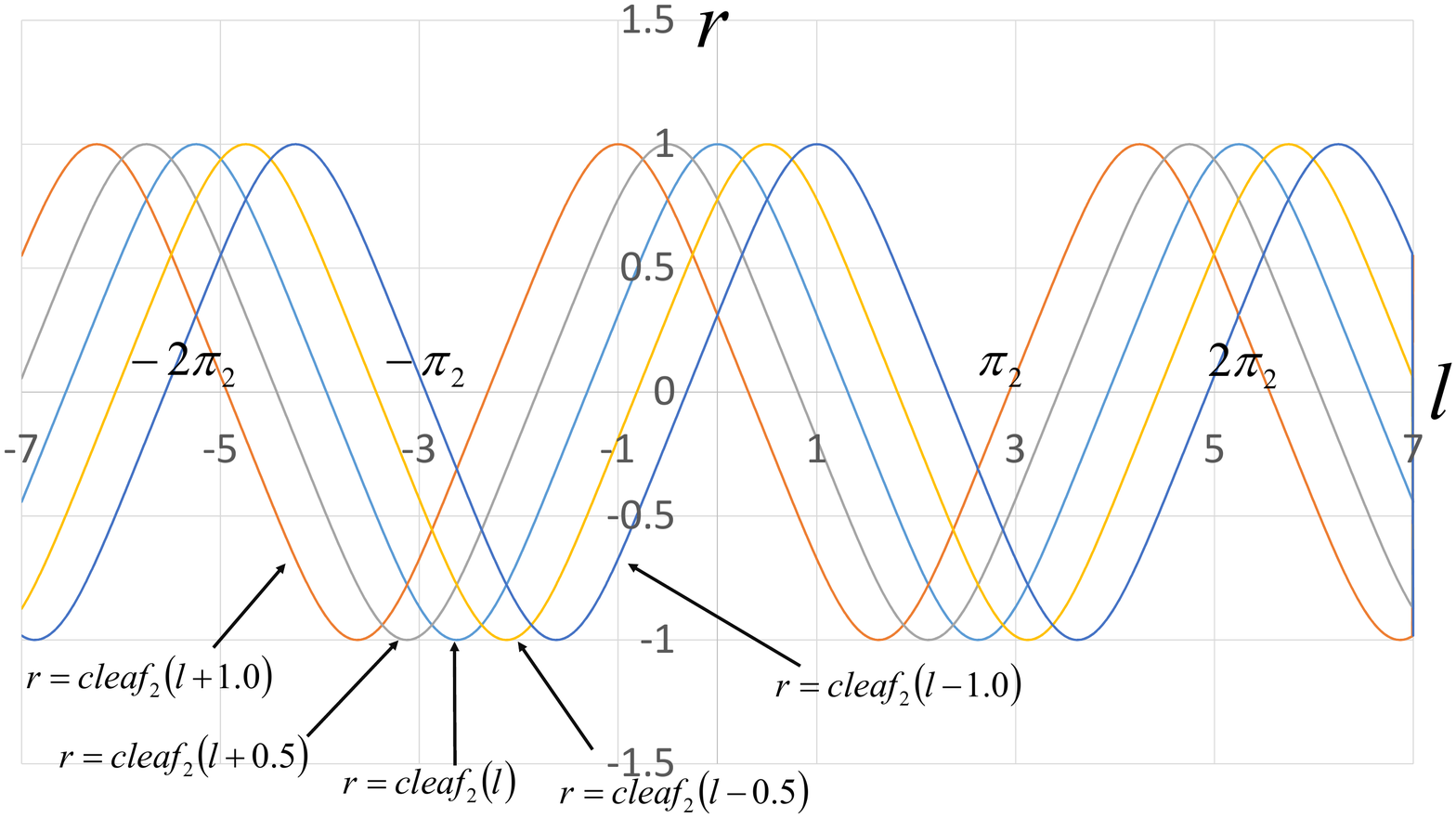}
\caption{ Translation of the curves of the function $\mathrm{cleaf}_2(l)$ obtained using the addition formulas with the basis $n = 2$ }
\label{fig4.1.5}          % Give a unique label
\end{center}
\end{figure*}

\begin{figure*}[tb]
\begin{center}
\includegraphics[width=0.75 \textwidth]{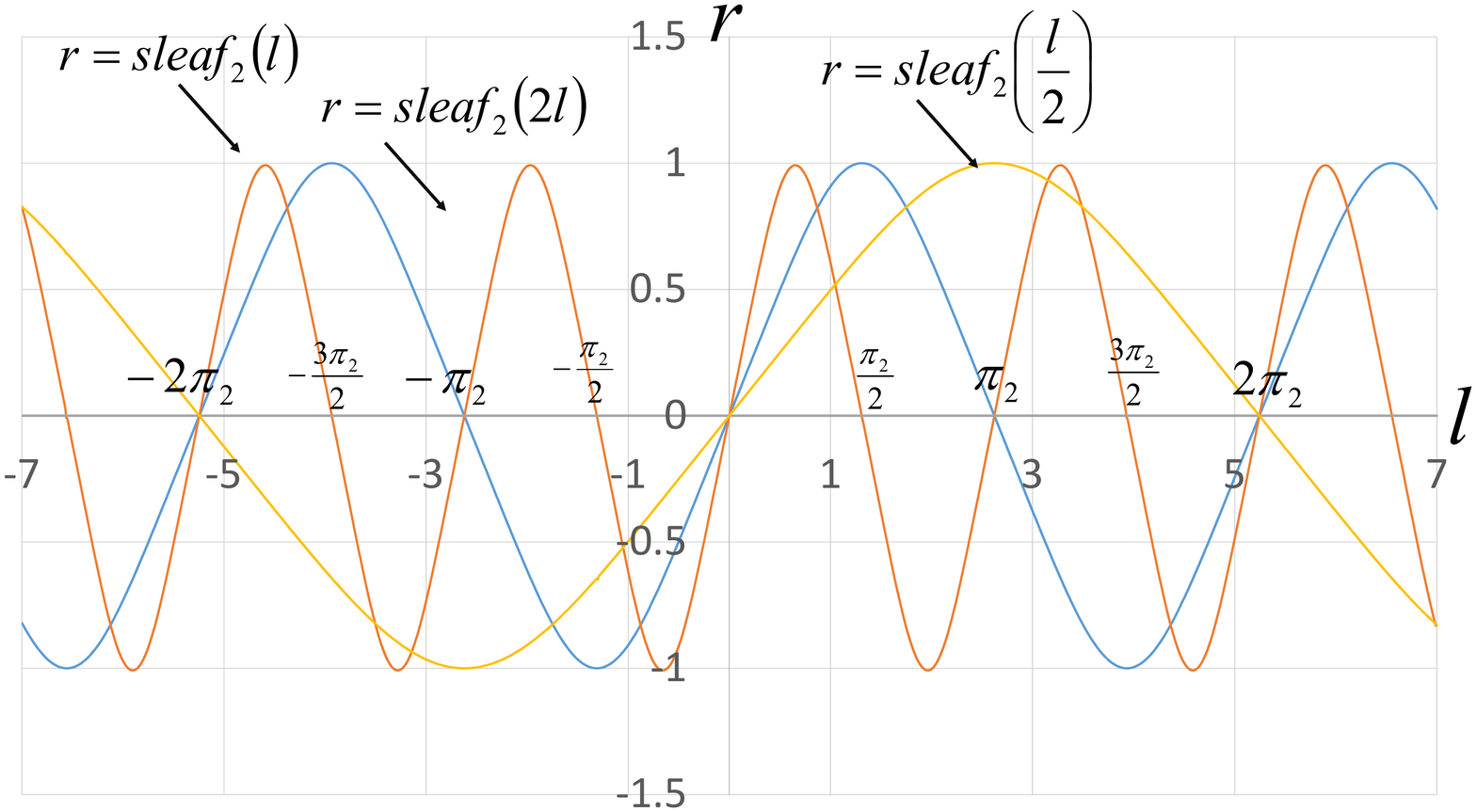}
\caption{ Translation of the curves of the functions $\mathrm{sleaf}_2(l)$, $\mathrm{sleaf}_2(2l)$, and $\mathrm{sleaf}_2(l/2)$  obtained using the addition formulas based on the basis $n = 2$ }
\label{fig4.1.6}        % Give a unique label
\end{center}
\end{figure*}

\begin{figure*}[tb]
\begin{center}
\includegraphics[width=0.75 \textwidth]{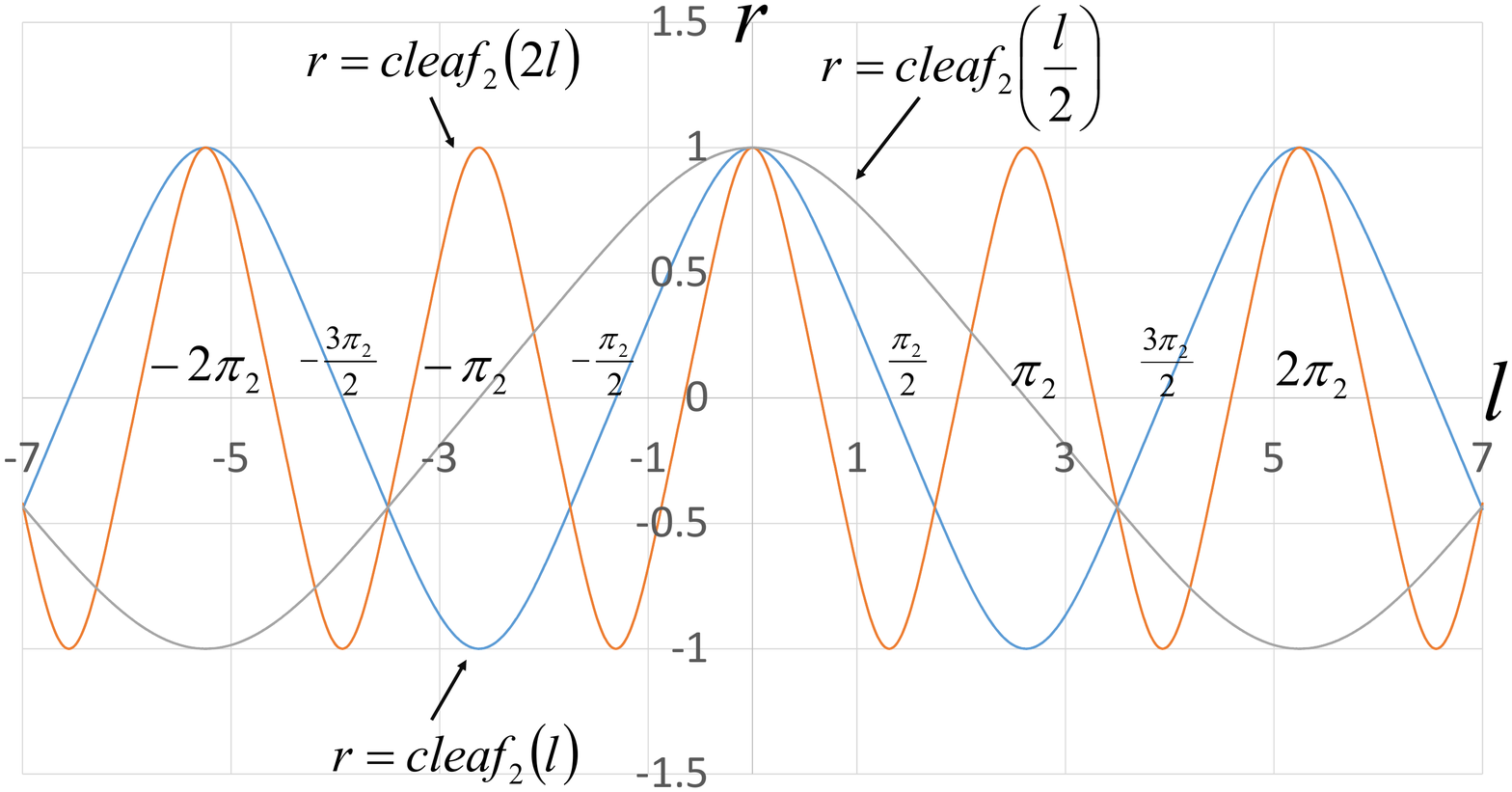}
\caption{ Translation of the curves of the functions $\mathrm{cleaf}_2(l)$, $\mathrm{cleaf}_2(2l)$, and $\mathrm{cleaf}_2(l/2)$  obtained using the addition formulas based on the basis $n = 2$ }
\label{fig4.1.7}   
\end{center}
\end{figure*}

The curves of the leaf functions $\mathrm{sleaf}_2(l)$ and $\mathrm{cleaf}_2 (l)$ are shown in Figs. \ref{fig4.1.4} and \ref{fig4.1.5}. Numerical data for these two leaf functions are summarized in Table 1. These curves are the same curves as those of the lemniscate elliptic functions $r=\mathrm{sl}(l)$ and $r=\mathrm{cl}(l)$. Using the addition formulas of Eq. (\ref{2.1.4}) $\sim$ Eq. (\ref{2.1.12}), the curves of the leaf functions $\mathrm{sleaf}_2(l)$ and $\mathrm{cleaf}_2(l)$ are translated in the direction of the axis $l$. Fig. \ref{fig4.1.6} shows graphs of the double-angle $\mathrm{sleaf}_2(2l)$ and the half-angle $\mathrm{sleaf}_2(l/2)$ obtained using Eqs. (\ref{3.1.4}) $\sim$ (\ref{3.1.5}) and Eqs. (\ref{3.2.3}) $\sim$ (\ref{3.2.4}) . Fig. \ref{fig4.1.7} shows graphs of the double-angle $\mathrm{cleaf}_2(2l)$ and the half-angle $\mathrm{cleaf}_2(l/2)$ obtained using Eqs. (\ref{3.1.6}) and (\ref{3.2.5}). The amplitude of the wave is 1 and one period of the function $\mathrm{cleaf}_2(l)$ is $2 \pi_2 (= 2 \times 2.622 \cdots)$. 

As shown in Fig. \ref{fig4.1.4}, curves $\mathrm{sleaf}_2(l)$ are translated using only the addition theorem, so that the period remains constant at $2 \pi_2$. On the contrary, as shown in Fig. \ref{fig4.1.6}, the period changes to $\pi_2$ and $4 \pi_2$, when the phase becomes $2l$ and $l / 2$, respectively. Similarly, as shown in Fig. \ref{fig4.1.5}, curves $\mathrm{cleaf}_2(l)$ are translated using only the addition theorem, so that the period remains constant at $2 \pi_2$. On the contrary, as shown in Fig. \ref{fig4.1.7}, the period changes to $\pi_2$ and $4 \pi_2$, when the phase becomes $2l$ and $l / 2$, respectively. Additionally, the leaf function can be expressed as the following trigonometric function:

\begin{equation}
\mathrm{sleaf}_n \Bigl(l+\frac{\pi_n}{2} \Bigr)=\mathrm{cleaf}_n (l) \: (n=1,2,3, \cdots)
\label{4.1.1}
\end{equation}

Using the Eq. (\ref{4.1.1}) and with constant  $\frac{\pi_2}{2}$, the waves are translated in the direction $l$. The curve shown in the Fig. \ref{fig4.1.6} represents the wave translated in the positive direction $l$, as shown in the Fig. \ref{fig4.1.7}. Similarly, the curve shown in the Fig. \ref{fig4.1.4} represents the wave translated in the positive direction $l$, as shown in the Fig. \ref{fig4.1.5}.

\begin{figure*}[tb]
\begin{center}
\includegraphics[width=0.75 \textwidth]{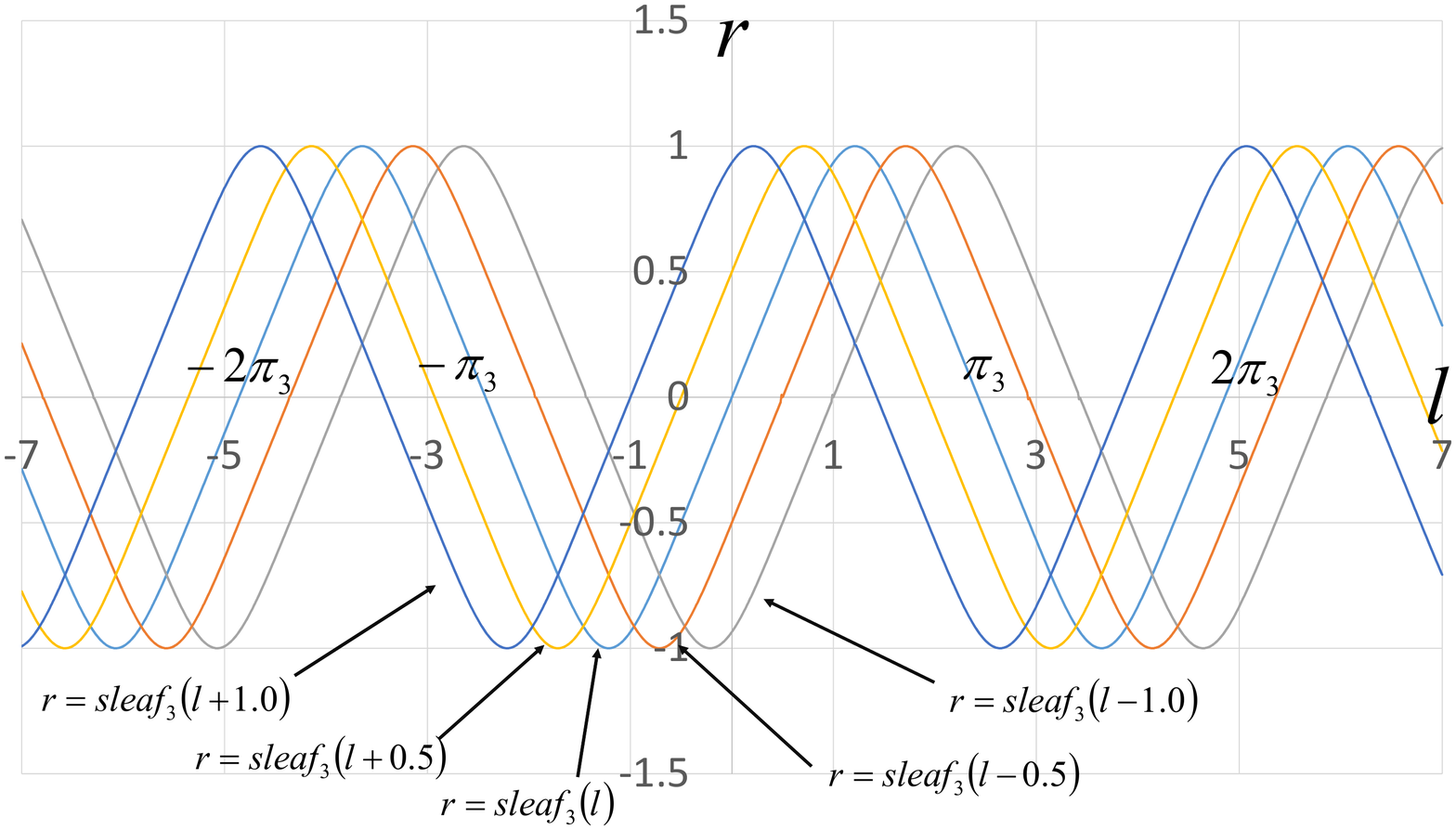}
\caption{ Translation of the curves of the function $\mathrm{sleaf}_3(l)$ obtained using the addition formulas with the basis $n = 3$ }
\label{fig4.1.8}   
\end{center}
\end{figure*}

\begin{figure*}[tb]
\begin{center}
\includegraphics[width=0.75 \textwidth]{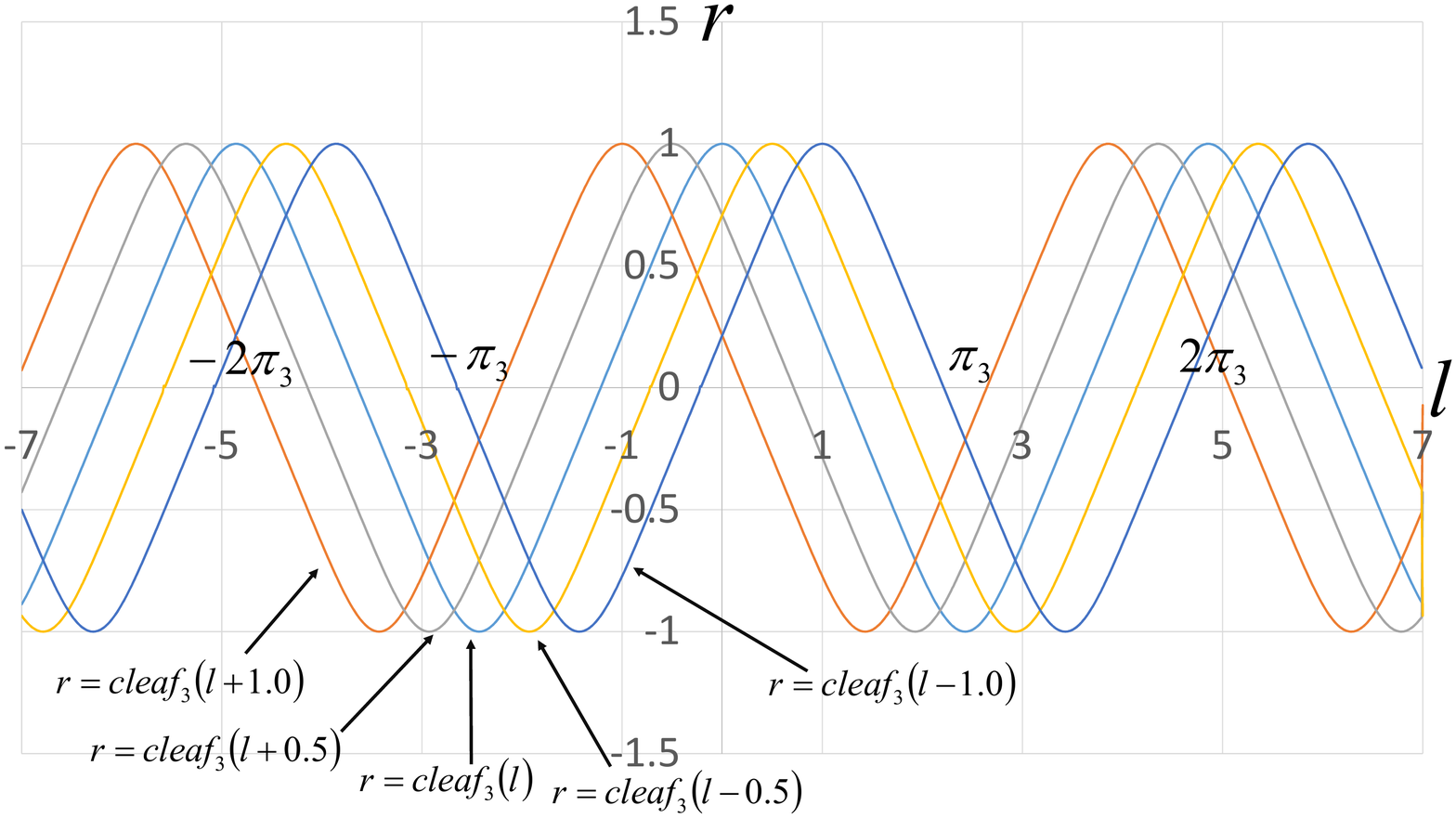}
\caption{ Translation of the curves of the function $\mathrm{cleaf}_3(l)$ obtained using the addition formulas with the basis $n = 3$ }
\label{fig4.1.9}   
\end{center}
\end{figure*}

\begin{figure*}[tb]
\begin{center}
\includegraphics[width=0.75 \textwidth]{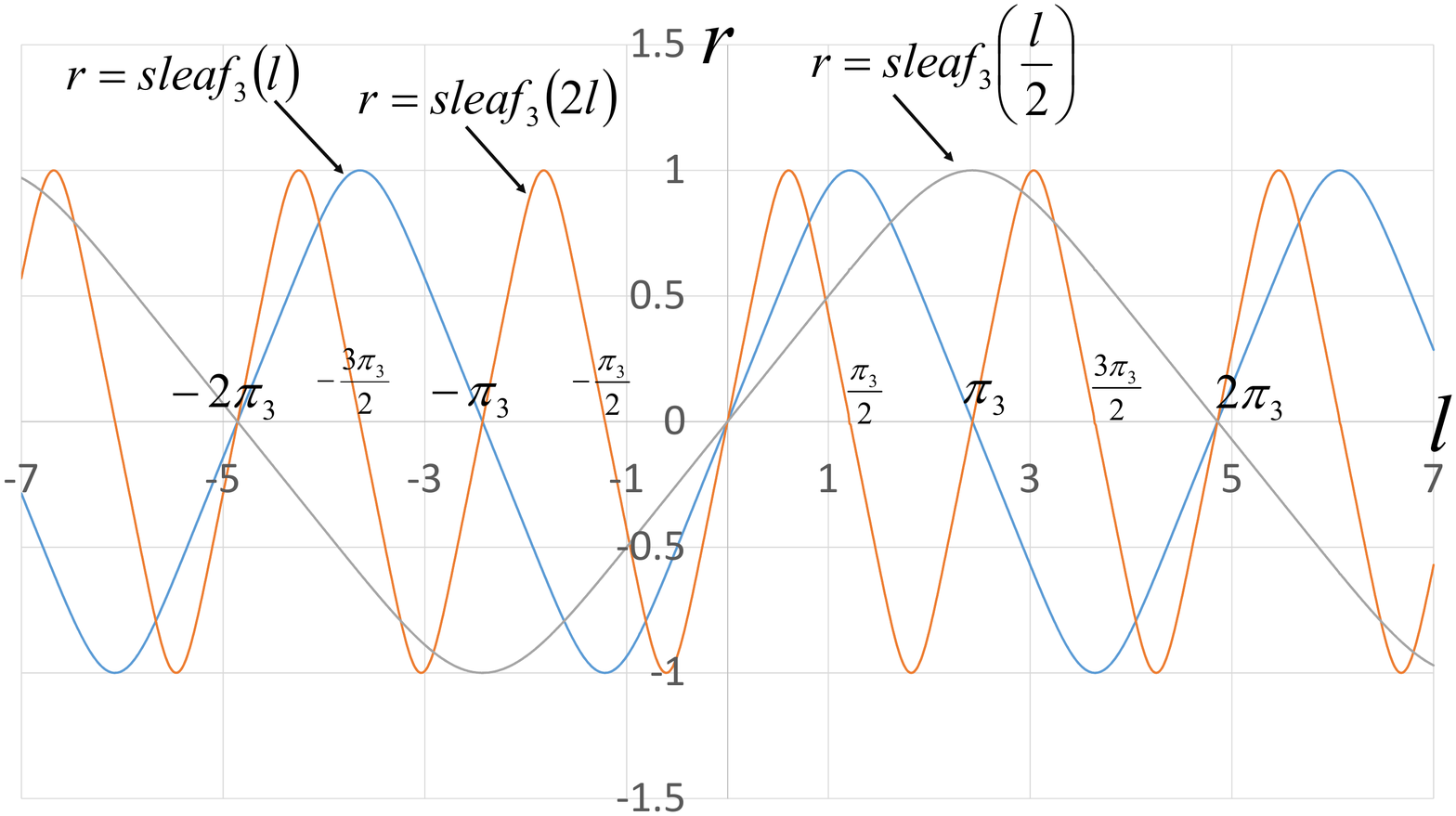}
\caption{ Translation of the curves of the functions $\mathrm{sleaf}_3(l)$, $\mathrm{sleaf}_3(2l)$, and $\mathrm{sleaf}_3(l/2)$  obtained using the addition formulas based on the basis $n = 3$ }
\label{fig4.1.10}   
\end{center}
\end{figure*}

\begin{figure*}[tb]
\begin{center}
\includegraphics[width=0.75 \textwidth]{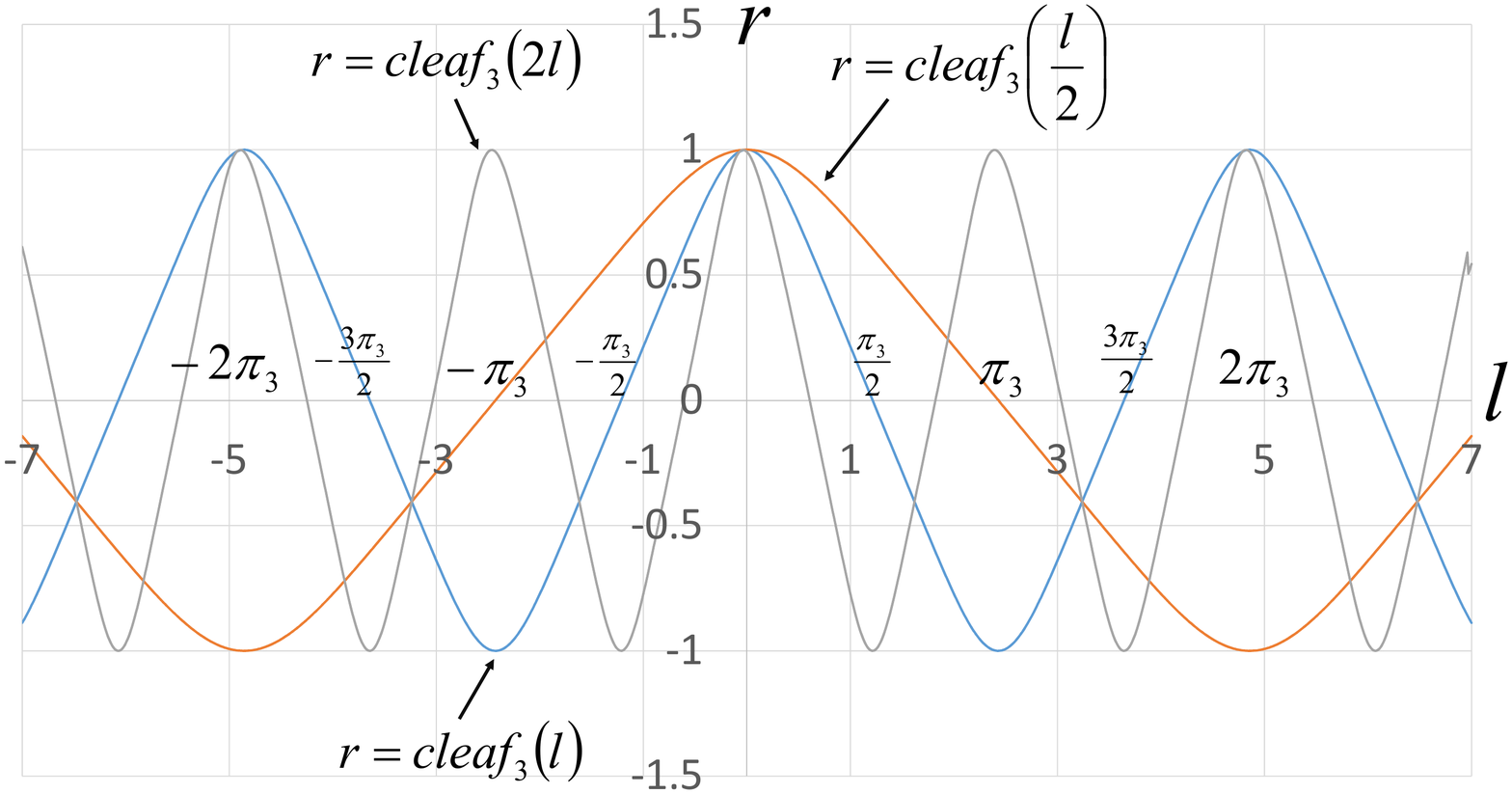}
\caption{ Translation of the curves of the functions $\mathrm{cleaf}_3(l)$, $\mathrm{cleaf}_3(2l)$, and $\mathrm{cleaf}_3(l/2)$  obtained using the addition formulas with the basis $n = 3$ }
\label{fig4.1.11}   
\end{center}
\end{figure*}

Next, the graph of the leaf function with the basis $n = 3$ is shown. The curves of the leaf functions $\mathrm{sleaf}_3(l)$ and $\mathrm{cleaf}_3(l)$ are shown in Figs. \ref{fig4.1.8} and \ref{fig4.1.9}. The horizontal and vertical axes represent the variables $l$ and $r$, respectively. The numerical data of the leaf functions $\mathrm{sleaf}_3(l)$ and $\mathrm{cleaf}_3(l)$ are summarized in Table \ref{tab:1}. The curves of the leaf functions $\mathrm{sleaf}_3(l)$ and $\mathrm{cleaf}_3(l)$ are translated in the direction of the axis $l$. These curves of the leaf functions were obtained using the addition formulas of Eqs. (\ref{2.1.14}) - (\ref{2.1.15}) and Eqs. (\ref{2.1.17}) - (\ref{2.1.18}). Fig. \ref{fig4.1.10} shows graphs of the double-angle $\mathrm{sleaf}_3(2l)$ and the half-angle $\mathrm{sleaf}_3(l/2)$ obtained using Eqs. (\ref{3.1.8}) $\sim$ (\ref{3.1.9}) and  Eqs. (\ref{3.2.6}) $\sim$ (\ref{3.2.7}).
Fig. \ref{fig4.1.11} shows graphs of the double-angle $\mathrm{cleaf}_3(2l)$ and the half-angle $\mathrm{cleaf}_3(l/2)$ obtained using Eq. (\ref{3.1.10}) and (\ref{3.2.8}). The amplitude of the wave is 1 and one period of the function $\mathrm{cleaf}_3(l)$ is $2\pi_3 (= 2 \times 2.429 \cdots)$.

  When the phase is doubled, the period is halved, and vice-versa. Even with a change in the phase, the amplitude remains constant at $1$, and initial condition  $\mathrm{cleaf}_3(0)=1$ is maintained at $l = 0$, as confirmed from the graph.

\subsection{Numerical Analysis of Hyperbolic Leaf Functions}
\label{Numerical Analysis of Hyperbolic Leaf Functions}
\begin{figure*}[tb]
\begin{center}
\includegraphics[width=0.75 \textwidth]{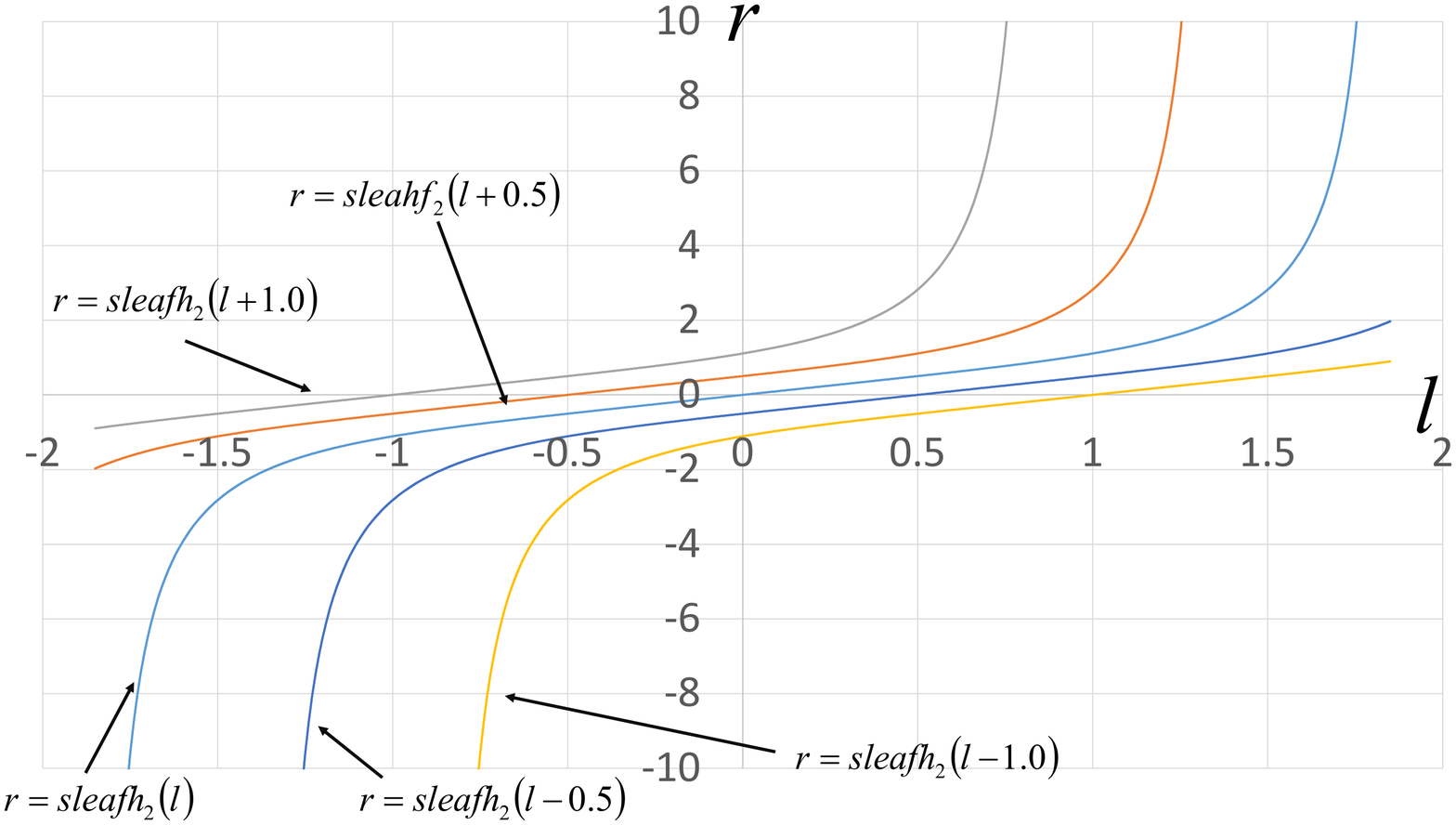}
\caption{ Translation of the curves of the function $\mathrm{sleafh}_2(l)$ obtained using the addition formulas with the basis $n =2$ }
\label{fig4.2.12}   
\end{center}
\end{figure*}

\begin{figure*}[tb]
\begin{center}
\includegraphics[width=0.75 \textwidth]{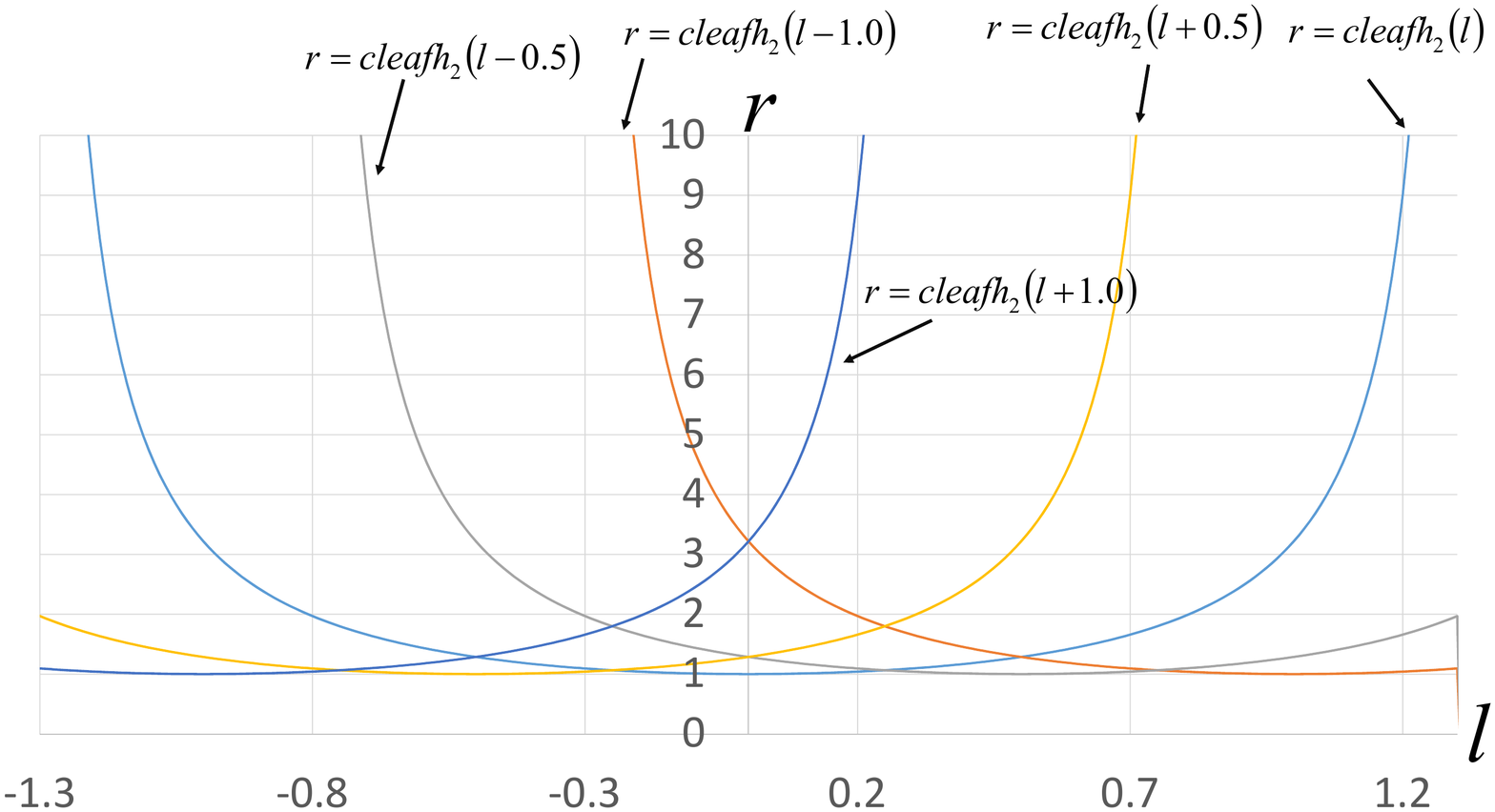}
\caption{ Translation of the curves of the function $\mathrm{cleafh}_2(l)$ obtained using the addition formulas with the basis $n =2$ }
\label{fig4.2.13}   
\end{center}
\end{figure*}

\begin{figure*}[tb]
\begin{center}
\includegraphics[width=0.75 \textwidth]{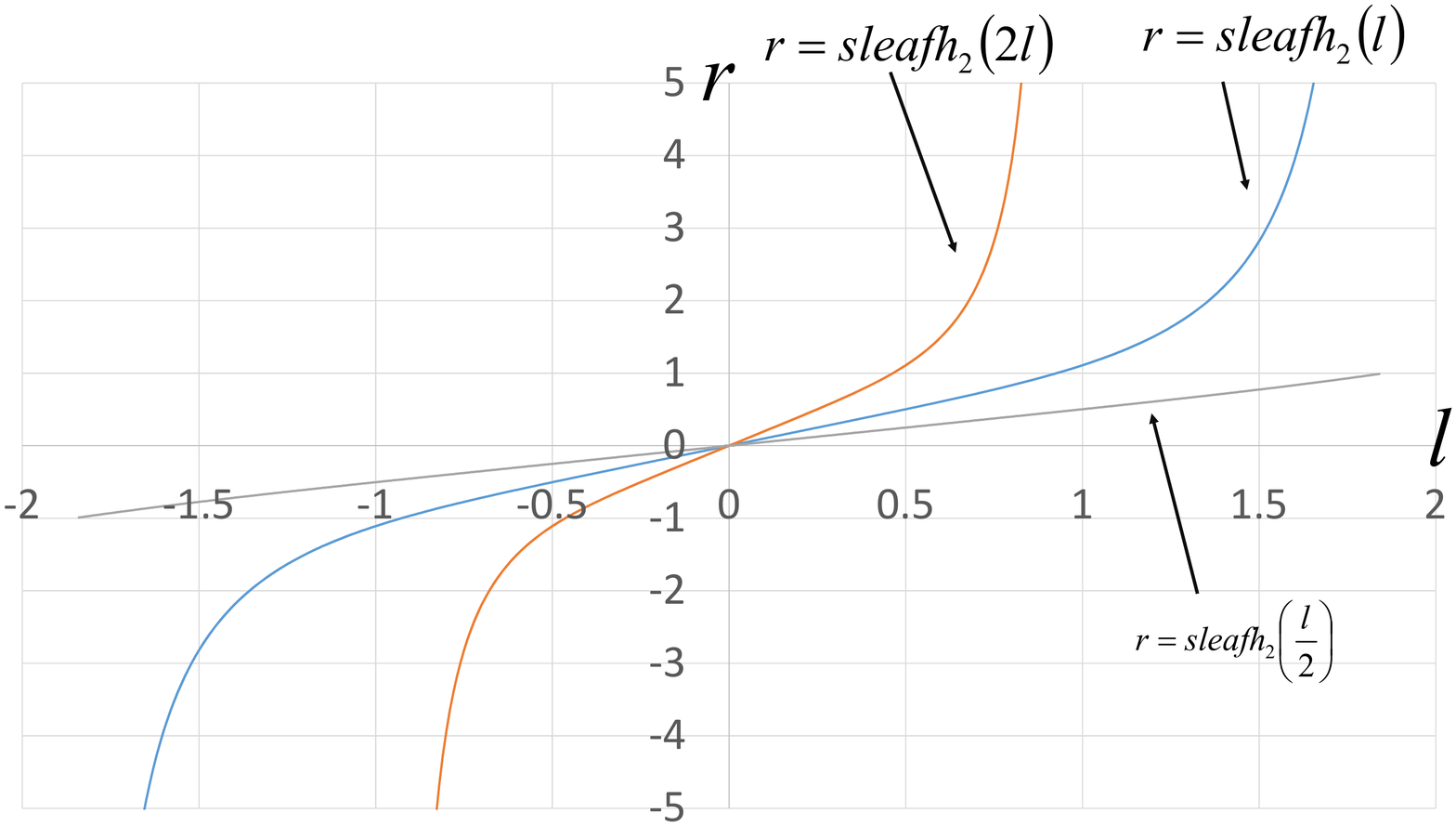}
\caption{ Translation of the curves of the functions $\mathrm{sleafh}_2(l)$, $\mathrm{sleafh}_2(2l)$ and $\mathrm{sleafh}_2(l/2)$  obtained using the addition formulas with the basis $n =2$ }
\label{fig4.2.14}   
\end{center}
\end{figure*}

\begin{figure*}[tb]
\begin{center}
\includegraphics[width=0.75 \textwidth]{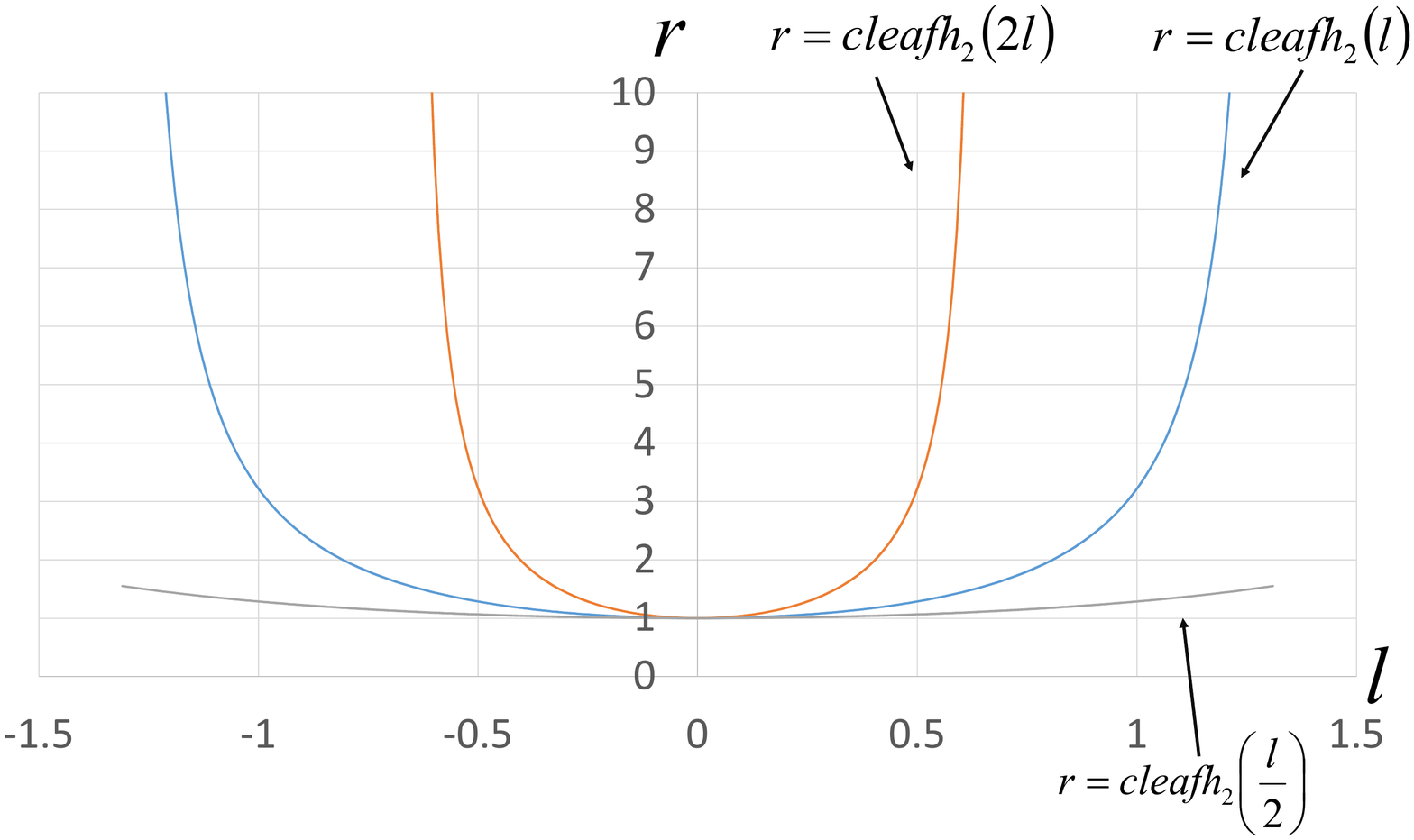}
\caption{  Translation of the curves of the functions $\mathrm{cleafh}_2(l)$, $\mathrm{cleafh}_2(2l)$, and $\mathrm{cleafh}_2(l/2)$  obtained using the addition formulas with the basis $n =2$ }
\label{fig4.2.15}   
\end{center}
\end{figure*}

\begin{figure*}[tb]
\begin{center}
\includegraphics[width=0.75 \textwidth]{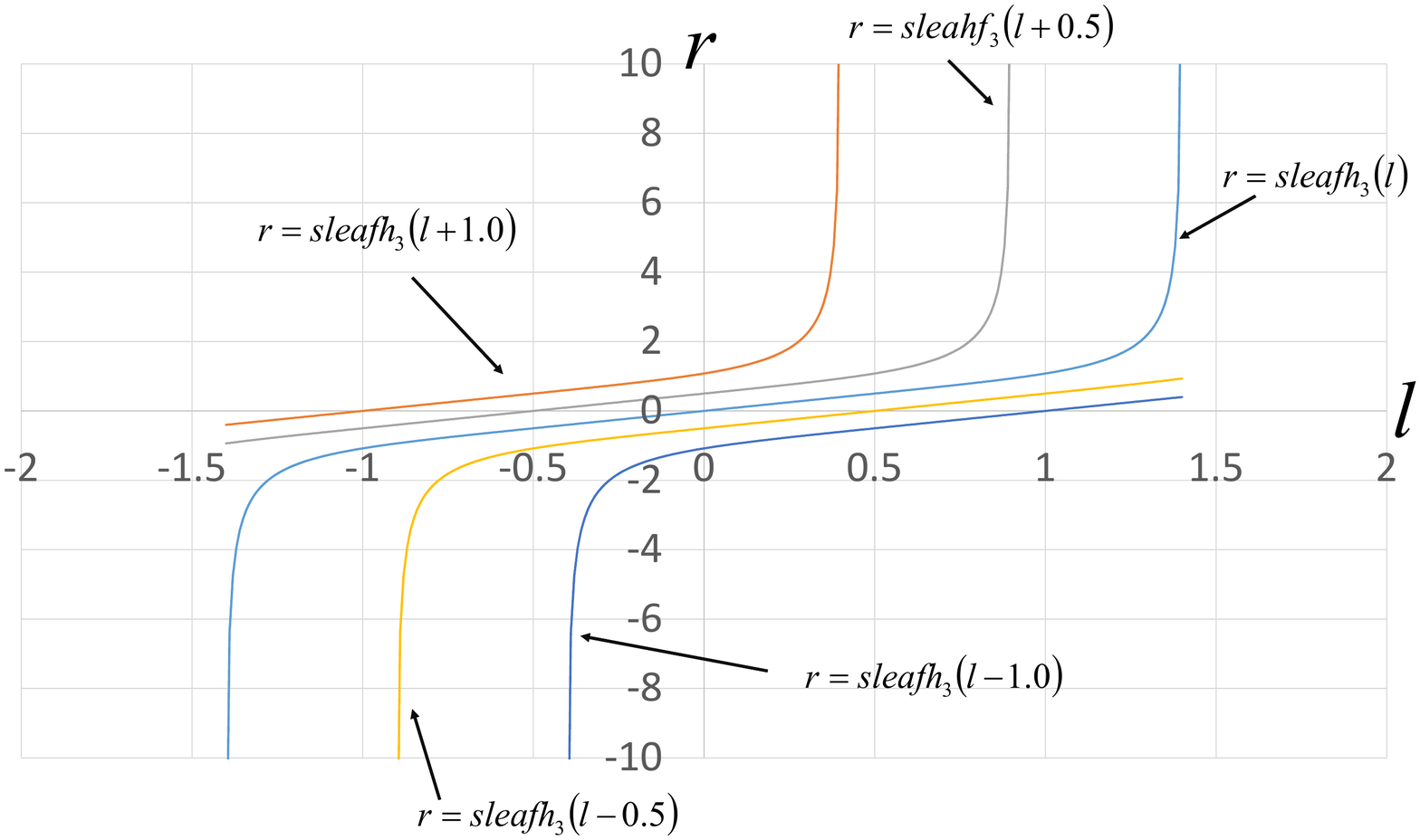}
\caption{ Translation of the curves of the function $\mathrm{sleafh}_3(l)$ obtained using the addition formulas with the basis $n =3$ }
\label{fig4.2.16}   
\end{center}
\end{figure*}

\begin{figure*}[tb]
\begin{center}
\includegraphics[width=0.75 \textwidth]{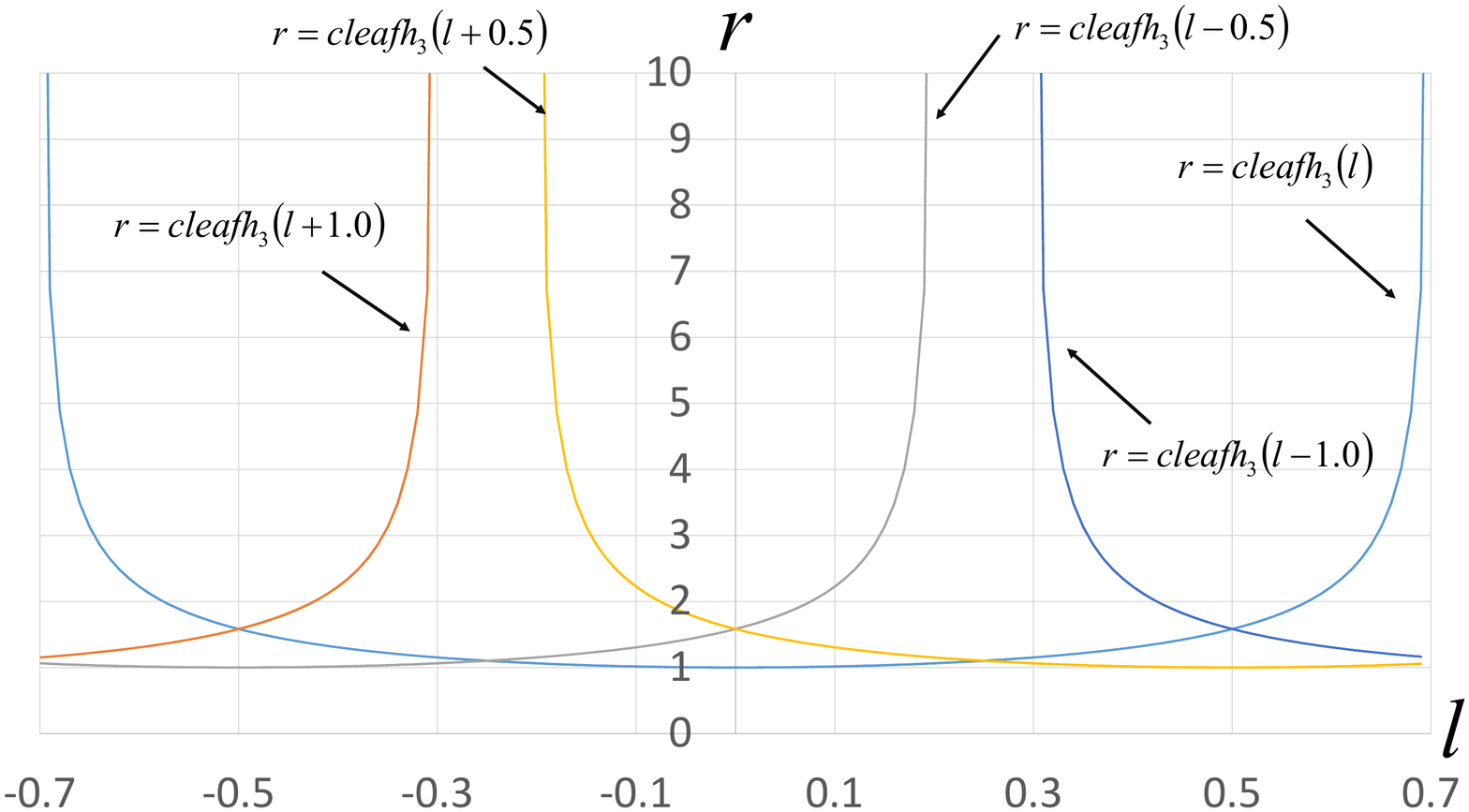}
\caption{ Translation of the curves of the functions $\mathrm{cleafh}_3(l)$ obtained using the addition formulas with the basis $n =3$ }
\label{fig4.2.17}   
\end{center}
\end{figure*}

\begin{figure*}[tb]
\begin{center}
\includegraphics[width=0.75 \textwidth]{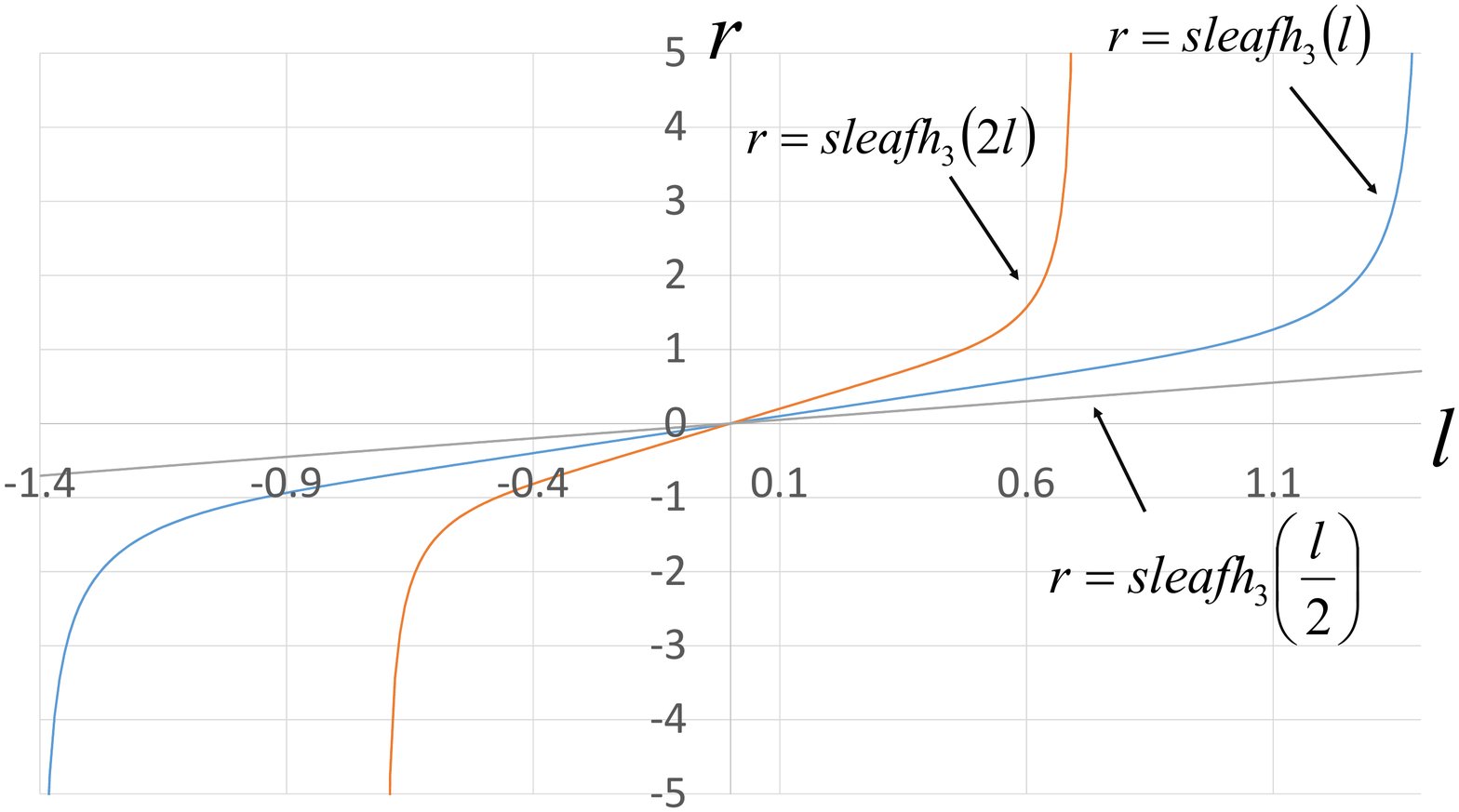}
\caption{ Translation of the curves of the functions $\mathrm{sleafh}_3(l)$, $\mathrm{sleafh}_3(2l)$, and $\mathrm{sleafh}_3(l/2)$  obtained using the addition formulas with the basis $n =3$ }
\label{fig4.2.18}   
\end{center}
\end{figure*}

The curves of the leaf functions $\mathrm{sleafh}_2(l)$ and $\mathrm{cleafh}_2(l)$ are shown in Figs. \ref{fig4.2.12} and \ref{fig4.2.13}. The horizontal and vertical axes represent the variables $l$ and $r$. The numerical data for the leaf functions $\mathrm{sleafh}_2(l)$ and $\mathrm{cleafh}_2(l)$  are summarized in Table \ref{tab:2}. Using the addition formulas of Eq. (\ref{2.2.10}) and the Eqs. (\ref{2.2.12}) $\sim$  (\ref{2.2.13}), the curves of the leaf functions $\mathrm{sleafh}_2(l)$ and $\mathrm{cleafh}_2(l)$ are translated in the direction $l$. Fig. \ref{fig4.2.14} shows graphs of the double-angle $\mathrm{sleafh}_2(2l)$ and the half-angle $\mathrm{sleafh}_2(l/2)$ obtained using Eq. (\ref{3.3.3}) and Eqs. (\ref{3.4.3}) $\sim$ (\ref{3.4.4}) . Fig. \ref{fig4.2.15} shows graphs of the double-angle $\mathrm{cleafh}_2(2l)$ and the half-angle $\mathrm{cleafh}_2(l/2)$ obtained using Eq. (\ref{3.3.4}) and Eqs. (\ref{3.4.5}) $\sim$ (\ref{3.4.6}). Limits exist for the functions $\mathrm{sleafh}_2(l)$ and $\mathrm{cleafh}_2(l)$, respectively (see Appendix F and Appendix G).
Next, curves of the leaf functions $\mathrm{sleafh}_3 (l)$ and $\mathrm{cleafh}_3(l)$ are shown in Figs. \ref{fig4.2.16} and \ref{fig4.2.17}. The horizontal and vertical axes represent the variables $l$ and $r$, respectively. The numerical data of the leaf functions $\mathrm{sleafh}_3(l)$ and $\mathrm{cleafh}_3(l)$ are summarized in Table 2. 
Using the addition formulas of Eq. (\ref{2.2.16}) and Eqs. (\ref{2.2.18}) $\sim$ (\ref{2.2.19}), the curves of the leaf functions $\mathrm{sleafh}_3(l)$ and $\mathrm{cleafh}_3(l)$ are translated in the direction $l$. Fig. \ref{fig4.2.18} shows graphs of the double-angle $\mathrm{sleafh}_3(2l)$ and the half-angle $\mathrm{sleafh}_3(l/2)$ obtained using Eq. (\ref{3.3.5}) and Eq. (\ref{3.4.7}). 

\begin{figure*}[tb]
\begin{center}
\includegraphics[width=0.75 \textwidth]{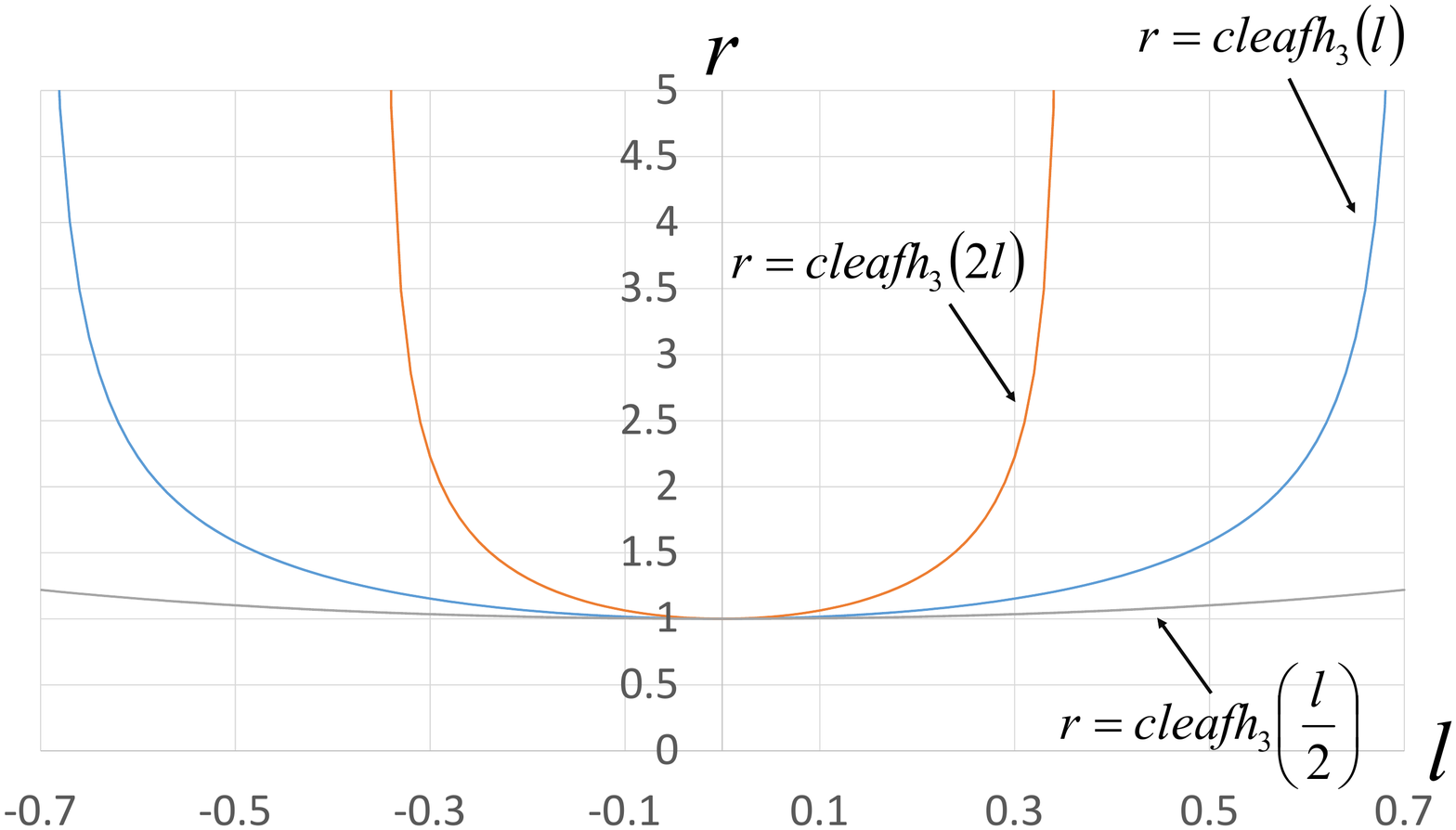}
\caption{ Translation of the curves of the functions $\mathrm{cleafh}_3(l)$, $\mathrm{cleafh}_3(2l)$, and $\mathrm{cleafh}_3(l/2)$  obtained using the addition formulas with the basis $n =3$ }
\label{fig4.2.19}   
\end{center}
\end{figure*}

Fig. \ref{fig4.2.19} shows graphs of the double-angle $\mathrm{cleafh}_3(2l)$ and the half-angle $\mathrm{cleafh}_3(l/2)$ obtained using Eqs (\ref{3.3.6}) and (\ref{3.4.8}). Limits exist in the functions $\mathrm{sleafh}_3(l)$ and $\mathrm{cleafh}_3(l)$, respectively. For the function $\mathrm{sleafh}_3(l)$, the limit exists at $\pm \zeta_3$ (see Appendix F for the constant $\zeta_3$). The curve of the function $\mathrm{sleafh}_3(l)$ monotonically increases in the domain $-\zeta_3 <l <\zeta_3$. In the case of the function $\mathrm{cleafh}_3(l)$, the limit exists at $\pm \eta_3$ (see Appendix G for the constant $\eta_3$). The domain of the function $\mathrm{cleafh}_3(l)$ is $-\eta_3 <l <\eta_3$.

\section{Conclusion}
\label{Conclusion}
Based on the analogy between the trigonometric and hyperbolic function, the hyperbolic leaf function paired with the leaf function was defined. The main conclusions can be summarized as follows:

$\cdot$ The relation equations between the leaf function and the hyperbolic leaf function were derived using imaginary numbers.

$\cdot$ The addition formulas of the hyperbolic leaf function were derived by using addition formulas of the leaf function with the basis $n = 1, 2, 3$. 

$\cdot$ For both the leaf function and hyperbolic leaf function for the basis $n = 1, 2, 3$, half-angle and double-angle formulas were derived using addition formulas 

As a future research topic, we will investigate whether the periodicity of the hyperbolic leaf function exists. In the case where the basis $n=2$, a limit exists in the hyperbolic function. By appropriately setting the initial conditions, the addition formulas for $n=2$ can be applied in all domains over the limit. Although the periodicity of the hyperbolic leaf function $n=2$ is evident, questions remain concerning the periodicity of the hyperbolic leaf function with $n=3$. In the case where the basis is $n = 3$, a limit also exists for the hyperbolic leaf function. However, the addition formulas of the hyperbolic leaf function cannot be applied outside of its domain. At basis $n = 3$, the periodicity of the hyperbolic leaf function is not observed.
 Another unaddressed issue is that the addition formulas of the leaf function with the basis $n=4$ or more are not known.  

\begin{table}
\centering
\caption{ Numerical data of the leaf functions }
\label{tab:1}   
\begin{tabular}{lllllll}
\hline\noalign{\smallskip}
$l$ & $\mathrm{sleaf}_1(l)$   & $\mathrm{cleaf}_1(l)$    & $\mathrm{sleaf}_2(l)$     & $\mathrm{cleaf}_2(l)$   & $\mathrm{sleaf}_3(l)$    & $\mathrm{cleaf}_3(l)$ \\
\noalign{\smallskip}\hline\noalign{\smallskip}
0.0	  &	0.000000000   &	1.000000000   &	0.000000000  	&	1.000000000  	&	0.000000000               	&	1.000000000	\\
0.1	&	0.099833417	&	0.995004165	&	0.099998987	&	0.990049602	&	0.099999991	&	0.98518434	\\
0.2	&	0.198669331	&	0.980066578	&	0.199967976	&	0.960781145	&	0.199999064	&	0.942809514	\\
0.3	&	0.295520207	&	0.955336489	&	0.299757126	&	0.913842132	&	0.299984331	&	0.878183695	\\
0.4	&	0.389418342	&	0.921060994	&	0.398978135	&	0.851676083	&	0.39988294	&	0.797825011	\\
0.5	&	0.479425539	&	0.877582562	&	0.496891146	&	0.777159391	&	0.499442694	&	0.70763201	\\
0.6	&	0.564642473	&	0.825335615	&	0.592307034	&	0.693234267	&	0.598009242	&	0.611978813	\\
0.7	&	0.644217687	&	0.764842187	&	0.683522566	&	0.602609146	&	0.694183101	&	0.513646507	\\
0.8	&	0.717356091	&	0.696706709	&	0.768312999	&	0.507563306	&	0.785387303	&	0.414175714	\\
0.9	&	0.78332691	&	0.621609968	&	0.844009686	&	0.409858439	&	0.867486256	&	0.314303714	\\
1.0	&	0.841470985	&	0.540302306	&	0.90768321	&	0.310738001	&	0.934767593	&	0.214323891	\\
1.1	&	0.89120736	&	0.453596121	&	0.956432623	&	0.210987025	&	0.980707849	&	0.114325366	\\
1.2	&	0.932039086	&	0.362357754	&	0.987748032	&	0.111027204	&	0.999692203	&	0.014325392	\\
1.3	&	0.963558185	&	0.267498829	&	0.999878378	&	0.011028912	&	0.989089542	&	-0.085674597	\\
1.4	&	0.98544973	&	0.169967143	&	0.99211532	&	-0.088970511	&	0.950392842	&	-0.185674048	\\
1.5	&	0.997494987	&	0.070737202	&	0.96491412	&	-0.188946955	&	0.888559535	&	-0.285663493	\\
1.6	&	0.999573603	&	-0.029199522	&	0.919815574	&	-0.288769649	&	0.810063642	&	-0.385583945	\\
1.7	&	0.99166481	&	-0.128844494	&	0.859192306	&	-0.388082304	&	0.720971617	&	-0.485219858	\\
1.8	&	0.973847631	&	-0.227202095	&	0.785891649	&	-0.486189025	&	0.6258955	&	-0.583992736	\\
1.9	&	0.946300088	&	-0.323289567	&	0.702864932	&	-0.581954203	&	0.527828311	&	-0.680635105	\\
2.0	&	0.909297427	&	-0.416146837	&	0.612857981	&	-0.673732946	&	0.428461029	&	-0.772765772	\\
2.1	&	0.863209367	&	-0.504846105	&	0.518203565	&	-0.759356014	&	0.328621294	&	-0.856486525	\\
2.2	&	0.808496404	&	-0.588501117	&	0.420721859	&	-0.836196738	&	0.228648563	&	-0.92628646	\\
2.3	&	0.745705212	&	-0.666276021	&	0.3217114	&	-0.90134206	&	0.128650882	&	-0.975673073	\\
2.4	&	0.675463181	&	-0.737393716	&	0.222003575	&	-0.951870972	&	0.028650956	&	-0.998769949	\\
2.5	&	0.598472144	&	-0.801143616	&	0.122054841	&	-0.985211764	&	-0.071349009	&	-0.992412076	\\
2.6	&	0.515501372	&	-0.856888753	&	0.022057545	&	-0.999513456	&	-0.171348665	&	-0.95749878	\\
2.7	&	0.42737988	&	-0.904072142	&	-0.077942171	&	-0.993943297	&	-0.2713412	&	-0.898594215	\\
2.8	&	0.33498815	&	-0.942222341	&	-0.177924624	&	-0.968828424	&	-0.371279371	&	-0.822087294	\\
2.9	&	0.239249329	&	-0.970958165	&	-0.277776677	&	-0.925599649	&	-0.470980082	&	-0.734191026	\\
3.0	&	0.141120008	&	-0.989992497	&	-0.37717265	&	-0.866554268	&	-0.569933963	&	-0.639752776	\\

\noalign{\smallskip}\hline
\end{tabular}
\end{table}

\begin{table}
\begin{center} 
\caption{ Numerical data of the hyperbolic leaf functions }
\label{tab:2}  
\begin{tabular}{rcccccc}
\hline\noalign{\smallskip}
$l$ & $\mathrm{sleafh}_1(l)$   & $\mathrm{cleafh}_1(l)$    & $\mathrm{sleafh}_2(l)$     & $\mathrm{cleafh}_2(l)$   & $\mathrm{sleafh}_3(l)$    & $\mathrm{cleafh}_3(l)$ \\
\noalign{\smallskip}\hline\noalign{\smallskip}
0.0	&	0.000000000	&	1.000000000	&	0.000000000	&	1.000000000	&	0.000000000	&	1.000000000	\\
0.1	&	0.10016675	&	1.005004168	&	0.100001013	&	1.010050409	&	0.100000009	&	1.015190873	\\
0.2	&	0.201336003	&	1.020066756	&	0.200032033	&	1.040819784	&	0.200000936	&	1.063219846	\\
0.3	&	0.304520293	&	1.045338514	&	0.300243205	&	1.094280966	&	0.300015671	&	1.152957367	\\
0.4	&	0.410752326	&	1.081072372	&	0.401026247	&	1.174155432	&	0.400117152	&	1.306327433	\\
0.5	&	0.521095305	&	1.127625965	&	0.503141445	&	1.286737533	&	0.500558986	&	1.583264962	\\
0.6	&	0.636653582	&	1.185465218	&	0.607861028	&	1.442514133	&	0.6020087	&	2.225120045	\\
0.7	&	0.758583702	&	1.255169006	&	0.717150413	&	1.659450947	&	0.705950043	&	21.4096535	\\
0.8	&	0.888105982	&	1.337434946	&	0.833926854	&	1.97019847	&	0.815368602	&$-$	\\
0.9	&	1.026516726	&	1.433086385	&	0.962467567	&	2.439868366	&	0.936017909	&$-$		\\
1.0	&	1.175201194	&	1.543080635	&	1.10910404	&	3.218148246	&	1.079143503	&$-$	\\
1.1	&	1.33564747	&	1.668518554	&	1.283479658	&	4.739635312	&	1.26866512	&$-$	\\
1.2	&	1.509461355	&	1.810655567	&	1.500980956	&	9.006830737	&	1.566095647	&$-$	\\
1.3	&	1.698382437	&	1.97091423	&	1.787828613	&	90.67397241	&	2.210887381	&$-$	\\
1.4	&	1.904301501	&	2.150898465	&	2.192926988	&$-$		&	15.13849028	&$-$ \\
\noalign{\smallskip}\hline
\end{tabular}
\end{center}
\end{table}

\begin{table}
\begin{center}
\caption{ Values of constants $\pi_n$ }
\label{tab:3}   
\begin{tabular}{ll}
\hline\noalign{\smallskip}
$n$ & $\pi_n$  \\
\noalign{\smallskip}\hline\noalign{\smallskip}
1	&	3.1415 $\cdots$	\\
2	&	2.6220 $\cdots$	\\
3	&	2.4286 $\cdots$	\\
$\cdots$	&	 $\cdots$	\\
\noalign{\smallskip}\hline
\end{tabular}
\end{center}
\end{table}

\begin{table}
\begin{center}
\caption{ Limits $\zeta_n$ of the hyperbolic leaf function $\mathrm{sleafh}_n(l)$ }
\label{tab:4}   
\begin{tabular}{ll}
\hline\noalign{\smallskip}
$n$ & $ \zeta_n$  \\
\noalign{\smallskip}\hline\noalign{\smallskip}
1	& Not applicable $\cdots$	\\
2	&	1.8540 $\cdots$	\\
3	&	1.4021 $\cdots$	\\
$\cdots$	&	 $\cdots$	\\
\noalign{\smallskip}\hline
\end{tabular}
\end{center}
\end{table}

\begin{table}
\begin{center}
\caption{ Limits $\eta_n$ of the hyperbolic leaf function $\mathrm{cleafh}_n(l)$ }
\label{tab:5}   
\begin{tabular}{ll}
\hline\noalign{\smallskip}
$n$ & $ \eta_n$  \\
\noalign{\smallskip}\hline\noalign{\smallskip}
1	& Not applicable $\cdots$	\\
2	&	1.31102 $\cdots$	\\
3	&	0.70109 $\cdots$	\\
$\cdots$	&	 $\cdots$	\\
\noalign{\smallskip}\hline
\end{tabular}
\end{center}
\end{table}

\appendix
\def\thesection{Appendix}
\section{A}
\label{Appendix A}
The relation equations with the basis $n=1$ are described. The relation equation between the leaf function $\mathrm{sleaf}_1(l)$ and the leaf function $\mathrm{cleaf}_1(l)$ is as follows:
\begin{equation}
 (\mathrm{sleaf}_1(l))^2+ (\mathrm{cleaf}_1(l))^2=1 \label{A1}
\end{equation}
The relation equation between the hyperbolic leaf function $\mathrm{sleafh}_1(l)$ and the hyperbolic leaf function $\mathrm{cleafh}_1(l)$ is as follows:
\begin{equation}
(\mathrm{cleafh}_1(l))^2 - (\mathrm{sleafh}_1(l))^2=1 \label{A2}
\end{equation}

\appendix
\def\thesection{Appendix}
\section{B}
\label{Appendix B}
The relation equations with the basis $n=2$ are described. The relation equation between the leaf function $\mathrm{sleaf}_2(l)$ and the leaf function $\mathrm{cleaf}_2(l)$ is as follows\cite{Kaz_cl}:
\begin{equation}
 (\mathrm{sleaf}_2(l))^2+ (\mathrm{cleaf}_2(l))^2+ (\mathrm{sleaf}_2(l))^2 \cdot (\mathrm{cleaf}_2(l))^2=1 \label{B1}
\end{equation}
The relation equation between the hyperbolic leaf function $\mathrm{sleafh}_2(l)$ and the hyperbolic leaf function $\mathrm{cleafh}_2(l)$ is as follows\cite{Kaz_sh, Kaz_ch}:
\begin{equation}
\mathrm{cleafh}_2(\sqrt{2} l) 
=\frac{1+ (\mathrm{sleafh}_2(l))^2}{1- (\mathrm{sleafh}_2(l))^2} \label{B2}
\end{equation}
The relation equation between the hyperbolic leaf function $\mathrm{cleaf}_2(l)$ and the hyperbolic leaf function $\mathrm{cleafh}_2(l)$ is as follows:
\begin{equation}
\mathrm{cleaf}_2(l) \cdot \mathrm{cleafh}_2(l)=1  \label{B3}
\end{equation}
The relation equation between the hyperbolic leaf function $\mathrm{sleaf}_2(l)$ and the hyperbolic leaf function $\mathrm{sleafh}_2(l)$ is as follows:
\begin{equation}
(\mathrm{sleaf}_2(\sqrt{2} l))^2 = \frac{ 2(\mathrm{sleafh}_2(l))^2 }{ 1+(\mathrm{sleafh}_2(l))^4 }  \label{B4}
\end{equation}

\appendix
\def\thesection{Appendix}
\section{C}
\label{Appendix C}
The relation equations with the basis $n=3$ are described. The relation equation between the leaf functions $\mathrm{sleaf}_3(l)$ and $\mathrm{cleaf}_3(l)$ is as follows \cite{Kaz_cl}:
\begin{equation}
 (\mathrm{sleaf}_3(l))^2+ (\mathrm{cleaf}_3(l))^2+ 2 (\mathrm{sleaf}_3(l))^2 \cdot (\mathrm{cleaf}_3(l))^2=1 \label{C1}
\end{equation}
The relation equation between the hyperbolic leaf functions $\mathrm{sleafh}_3(l)$ and  $\mathrm{cleafh}_3(l)$ is as follows\cite{Kaz_sh, Kaz_ch}:
\begin{equation}
  (\mathrm{cleafh}_3(l))^2-(\mathrm{sleafh}_3(l))^2- 2 (\mathrm{sleafh}_3(l))^2 \cdot (\mathrm{cleafh}_3(l))^2=1 \label{C2}
\end{equation}

\appendix
\def\thesection{Appendix}
\section{D}
\label{Appendix D}
Using the imaginary number, the relations between the leaf function and hyperbolic leaf function are described in the works\cite{Kaz_sh, Kaz_ch}.
To derive the relation between these two functions, the following equation is defined:
\begin{equation}
r=i \cdot u \label{D1}
\end{equation}
The symbol $i$ represents the imaginary number. Substituting the preceding equation  yields the following:
\begin{equation}
l=\int_{0}^{i \cdot u} \frac{\mathrm{d}t}{\sqrt{1-t^{2n}}}
(=\mathrm{arcsleaf}_n(i \cdot u)) \label{D2}
\end{equation}
Here, the parameter $t$ is replaced with $i \cdot \xi$ ($t=i \cdot \xi$). In the case where $t=0$, $\xi$ is zero. In the case where $t=i \cdot u $, $\xi$ is $u$. Thus, the following equation is obtained:
\begin{equation}
l=\int_{0}^{u} \frac{i \cdot \mathrm{d} \xi }{\sqrt{1-( i \cdot \xi )^{2n}}}
=i \cdot \int_{0}^{u} \frac{ \mathrm{d} \xi }{\sqrt{ 1- i^{2n} \cdot \xi^{2n} }} \label{D3}
\end{equation}
Let $n$ be an odd number, that is, $n = 2m-1 (m =1, 2, 3, \cdots)$. The following equation is then obtained:
\begin{equation}
l=i \cdot \int_{0}^{u} \frac{\mathrm{d} \xi }{\sqrt{ 1- i^{2n} \cdot \xi^{2n} }}
=i \cdot \int_{0}^{u} \frac{\mathrm{d} \xi }{\sqrt{ 1 + \xi^{2n} }}
=i \cdot \mathrm{asleafh}_n(u) \label{D4}
\end{equation}
The following equation is obtained based on the preceding equation as follows:
\begin{equation}
\mathrm{sleafh}_n \Bigl(\frac{l}{i} \Bigr)=u \label{D5}
\end{equation}
\begin{equation}
\mathrm{sleafh}_n (-i \cdot l )=u \label{D6}
\end{equation}
Here, the leaf function $\mathrm{sleafh}_n (l)$ has the following relation\cite{Kaz_sh}:
\begin{equation}
\mathrm{sleafh}_n(- l)=-\mathrm{sleafh}_n(l) \label{D7}
\end{equation}
Eq. (\ref{D5}) can be transformed as follows:
\begin{equation}
-\mathrm{sleafh}_n(i \cdot l)=u \label{D8}
\end{equation}
The following equation is obtained using Eq. (\ref{D2}) and Eq. (\ref{D8}):
\begin{equation}
\mathrm{sleaf}_n(l)=-i \cdot \mathrm{sleafh}_n(i \cdot l) \label{D9}
\end{equation}
Next, let us consider the case where $n$ is an even number. In the case where $n=2m (m =1, 2, 3 \cdots)$, the following equation is obtained:
\begin{equation}
l=i \cdot \int_{0}^{u} \frac{\mathrm{d} \xi }{\sqrt{ 1- i^{2n} \cdot \xi^{2n} }}
=i \cdot \int_{0}^{u} \frac{\mathrm{d} \xi }{\sqrt{ 1 - \xi^{2n} }}
=i \cdot \mathrm{arcsleaf}_n(u) \label{D10}
\end{equation}
The following equation is obtained: 
\begin{equation}
\mathrm{sleaf}_n \Bigl(\frac{l}{i} \Bigr)=u \label{D11}
\end{equation}
\begin{equation}
\mathrm{sleaf}_n (-i \cdot l)=u \label{D12}
\end{equation}
Here, the leaf function $\mathrm{sleaf}_n(l)$ has the following relation\cite{Kaz_sl}:
\begin{equation}
\mathrm{sleaf}_n(- l)=-\mathrm{sleaf}_n(l) \label{D13}
\end{equation}
Eq. (\ref{D12}) can be expressed as follows:
\begin{equation}
-\mathrm{sleaf}_n(i \cdot l)=u \label{D14}
\end{equation}
The following equation is obtained using Eq. (\ref{D2}) and Eq. (\ref{D14}):
\begin{equation}
\mathrm{sleaf}_n(l)=-i \cdot \mathrm{sleaf}_n(i \cdot l) \label{D15}
\end{equation}
In the case where $n$ is an even number, the following equation is also derived:
\begin{equation}
\mathrm{sleafh}_n(l)=-i \cdot \mathrm{sleafh}_n(i \cdot l) \label{D16}
\end{equation}
Next, the equation can be transformed as follows:
\begin{equation}
l= \int_{1}^{r} \frac{\mathrm{d}t }{\sqrt{ t^{2n} -1 } }
=  \int_{1}^{r} \frac{\mathrm{d}t }{i \sqrt{ 1- t^{2n} } }
=\frac{1}{i} \cdot \int_{1}^{r} \frac{\mathrm{d}t }{\sqrt{ 1- t^{2n} } }
=\frac{1}{i}  \mathrm{arccleaf}_n(r) \label{D17}
\end{equation}
The following equation is obtained by the Eq.  (\ref{D17}):
\begin{equation}
r=\mathrm{cleaf}_n(i \cdot l) \label{D18}
\end{equation}
The following equation is also obtained by the Eq.  (\ref{D17}): 
\begin{equation}
r=\mathrm{cleafh}_n(l) \label{D19}
\end{equation}
The following equation is obtained using Eq. (\ref{D18}) and Eq. (\ref{D19}):
\begin{equation}
\mathrm{cleaf}_n(i \cdot l)=\mathrm{cleafh}_n(l) \label{D20}
\end{equation}
Alternatively, the following equation is obtained by substituting $i \cdot l$ into $l$:
\begin{equation}
\mathrm{cleaf}_n(l)=\mathrm{cleafh}_n(i \cdot l) \label{D21}
\end{equation}
In the preceding equation, the following equation is applied:
\begin{equation}
\mathrm{cleaf}_n(l)=\mathrm{cleaf}_n(-l) \label{D22}
\end{equation}

\appendix
\def\thesection{Appendix}
\section{E}
\label{Appendix E}
The constants $\pi_n$ are defined as follows\cite{Kaz_sl, Kaz_cl}:
\begin{equation}
\pi_n=2 \int_{0}^{1} \frac{1}{ \sqrt{1-t^{2n}} } \mathrm{d}t (n=1,2,3 \cdots) \label{E1}
\end{equation}
In the case where $n=1$, the constant $\pi_1$ represents the circular constant $\pi$. The constants $\pi_n (n=1,2,3 \cdots)$ are summarized in Table \ref{tab:3}.

\appendix
\def\thesection{Appendix}
\section{F}
\label{Appendix F}
Except for the basis $n=1$, the limit of the variable $l$ exists in the hyperbolic leaf function $\mathrm{sleafh}_n(l)$\cite{Kaz_sh}. The limit with the basis $n$ is defined as $\zeta_n$. The limit $\zeta_n$ is obtained by the following equation:
\begin{equation}
\zeta_n= \int_{0}^{\infty} \frac{1}{ \sqrt{1+t^{2n}} }\mathrm{d}t (n=2,3 \cdots) \label{F1}
\end{equation}
The constants $\zeta_n (n=2,3 \cdots)$ are summarized in Table \ref{tab:4}.

\appendix
\def\thesection{Appendix}
\section{G}
\label{Appendix G}
Except for the basis $n=1$, the limit of the variable $l$ exists in the hyperbolic leaf function $\mathrm{cleafh}_n(l)$\cite{Kaz_ch}. The limit with the basis $n$ is defined as $\eta_n$. The limit $\eta_n$ is obtained using the following equation:
\begin{equation}
\eta_n= \int_{1}^{\infty} \frac{1}{ \sqrt{t^{2n}-1} }\mathrm{d}t  (n=2,3 \cdots) \label{G1}
\end{equation}
The constants $\eta_n (n=2,3 \cdots)$ are summarized in Table \ref{tab:5}.

\appendix
\def\thesection{Appendix}
\section{H}
\label{Appendix H}
\begin{figure*}[tb]
\includegraphics[width=0.75 \textwidth]{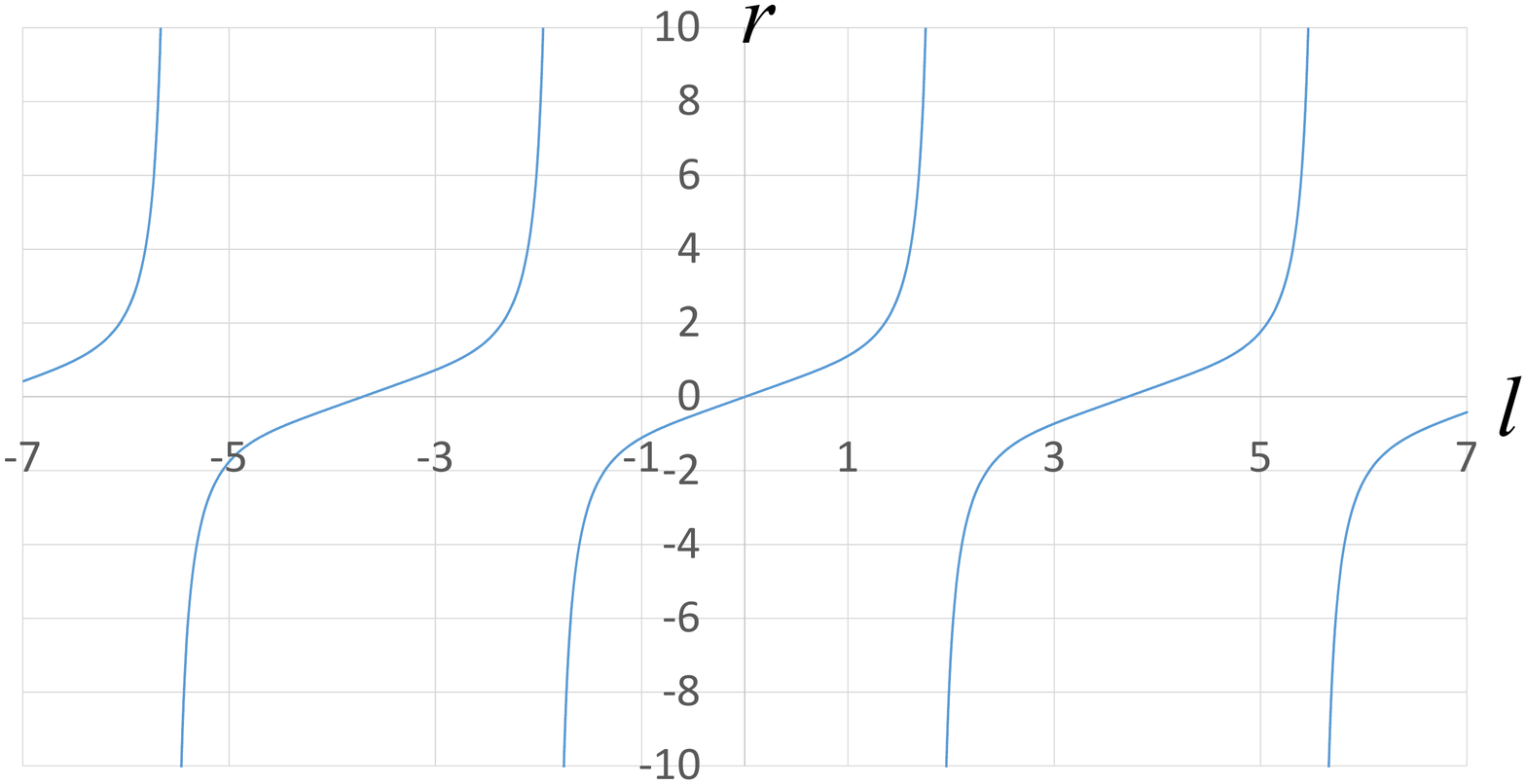}
\caption{Curves of the extended hyperbolic leaf function $\mathrm{sleafh}_2(l)$ for the initial conditions: Eq. (\ref{H1}) and Eq. (\ref{H2})}
\label{figH.1}       % Give a unique label
\end{figure*}

\begin{figure*}[tb]
\includegraphics[width=0.75 \textwidth]{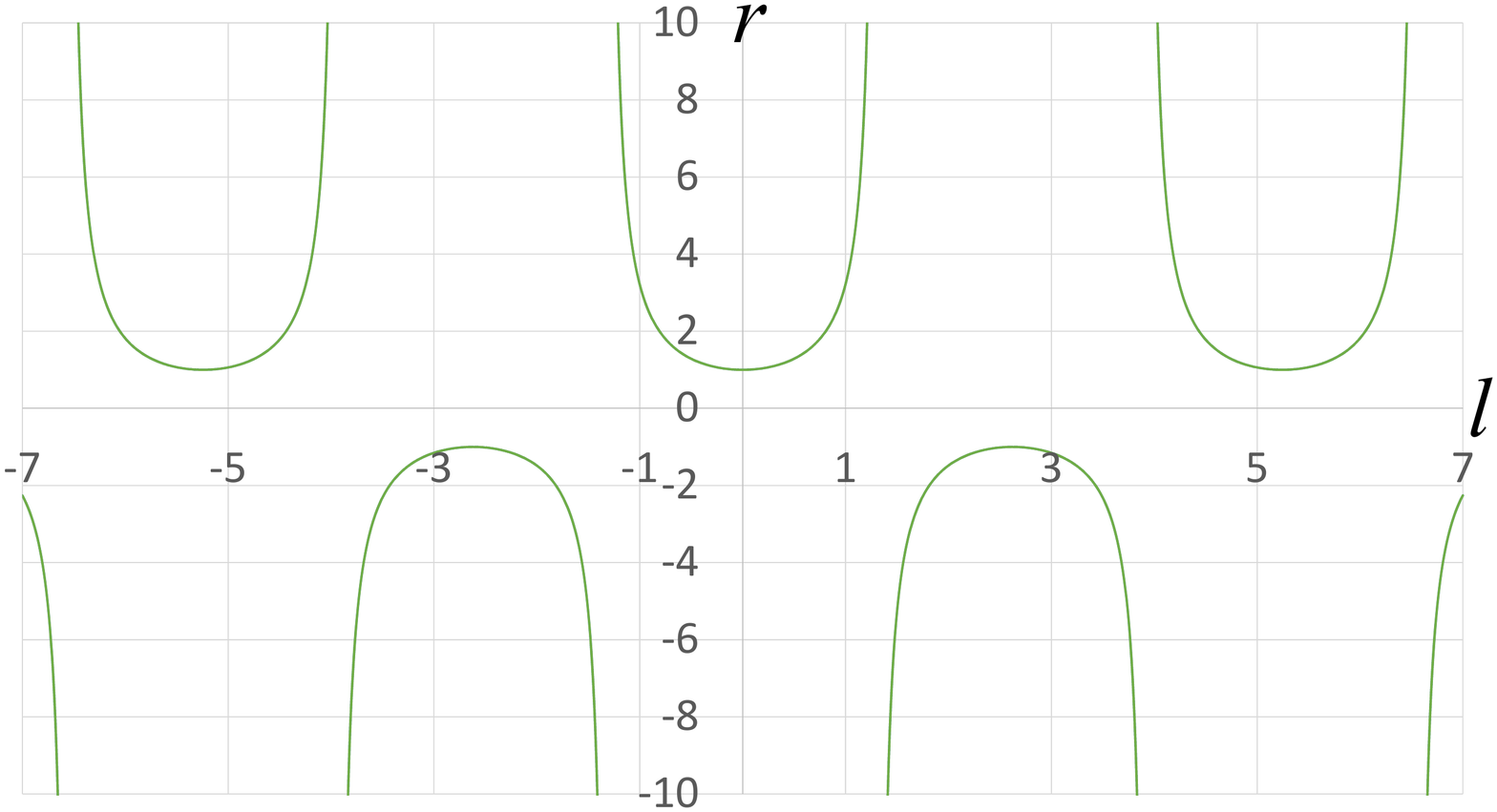}
\caption{Curves of the extended hyperbolic leaf function $\mathrm{cleafh}_2(l)$ for the initial conditions: Eq. (\ref{H3}), Eq. (\ref{H4}), and Eq. (\ref{H5})}
\label{figH.2}  
\end{figure*}

The function $\mathrm{sleafh}_n(l)$ and $\mathrm{cleafh}_n(l)$ have limits. The domains of the variable $l$ are defined as Eq. (\ref{1.1.8}) and Eq. (\ref{1.1.12}), respectively.  Therefore, the values of the hyperbolic leaf function cannot be defined under the domain $|l|>\zeta_n$ in the function $\mathrm{sleafh}_n(l)$ or $|l|>\eta_n$ in the function $\mathrm{cleafh}_n(l)$. In the case where $n=1$, the limits do not exist in the hyperbolic leaf function as $\mathrm{sleafh}_1(l)$ and $\mathrm{cleafh}_1(l)$ represent $\mathrm{sinh}(l)$ and $\mathrm{cosh}(l)$, respectively.
In the case where $n=2$ ($\mathrm{sleafh}_2(l)$ and $\mathrm{cleafh}_2(l)$), the initial values of the variables $r(0)$ and $dr(0)/dt$ are defined by Eqs. (\ref{1.1.9}) and (\ref{1.1.10}), or Eq. (\ref{1.1.13}) and (\ref{1.1.14}). The initial values in the function $\mathrm{sleafh}_2(l)$ are redefined as follows:
\begin{equation}
r(2m \zeta_2)=\mathrm{sleafh}_2(2m \zeta_2)=0 \label{H1}
\end{equation}
\begin{equation}
\frac{\mathrm{d}r(2m \zeta_2)}{\mathrm{d}l}=\frac{\mathrm{d}}{\mathrm{d}l} \mathrm{sleafh}_2(2m \zeta_2)=1  \label{H2}
\end{equation}
The initial values of the function $\mathrm{cleafh}_2(l)$ are redefined as follows:
\begin{equation}
r(4m \eta_2)=\mathrm{cleafh}_2(4m \eta_2)=1   \label{H3}
\end{equation}
\begin{equation}
r((4m-2) \eta_2)=\mathrm{cleafh}_2((4m-2)  \eta_2)=-1   \label{H4}
\end{equation}
\begin{equation}
\frac{\mathrm{d}r(2m \eta_2)}{\mathrm{d}l}=\frac{\mathrm{d}}{\mathrm{d}l} \mathrm{cleafh}_2(2m \eta_2)=0  \label{H5}
\end{equation}
The variable $m$ represents an integer. The graph based on these definitions is shown in Fig. \ref{figH.1} ($\mathrm{sleafh}_2(l)$ ) and Fig. \ref{figH.2} ($\mathrm{cleafh}_2(l)$ ), respectively. Such definitions are consistent for all the formulas such as the addition, double-angle, and half-angle formulas. These formulas work under all domains.
In the case $n=2$, the hyperbolic leaf functions can be extended for all domains. In the case where $n=3$ in the hyperbolic leaf function, the addition, double-angle, and half-angle formulas do not work in the domain $|l|>\zeta_3$ of $|l|>\eta_3$, even if the initial conditions are defined by equations such as $r(2m \zeta_3)=\mathrm{sleafh}_3(2m \zeta_3)=0$. In the case where $n \geqq 3$, the values of $\mathrm{sleafh}_n(l)$ and $\mathrm{cleafh}_n$ are unknown for the domain $|l|>\zeta_n$ of $|l|>\eta_n$. 

\appendix
\def\thesection{Appendix}
\section{I}
\label{Appendix I}
Eq. (\ref{2.1.13}) is set as follows:
\begin{equation}
\{ g(l_1,l_2) \}^2=\frac{ \{p_1(l_1,l_2) \}^2 }{p_3(l_1,l_2)}+\frac{ \{p_2(l_1,l_2) \}^2 }{p_3(l_1,l_2)}  \label{I1}
\end{equation}
\begin{equation}
g(l_1,l_2)= \mathrm{sleaf}_3(l_1+l_2)  \label{I2}
\end{equation}
\begin{equation}
p_1(l_1,l_2) = \mathrm{sleaf}_3(l_1) \frac{ \partial \mathrm{sleaf}_3(l_2) }{ \partial l_2}+ \mathrm{sleaf}_3(l_2) \frac{ \partial \mathrm{sleaf}_3(l_1) }{ \partial l_1}  \label{I3}
\end{equation}
\begin{equation}
p_2(l_1,l_2) = (\mathrm{sleaf}_3(l_1))^3  \mathrm{sleaf}_3(l_2) - \mathrm{sleaf}_3(l_1)  (\mathrm{sleaf}_3(l_2))^3   \label{I4}
\end{equation}
\begin{equation}
p_3(l_1,l_2) = 1+4 (\mathrm{sleaf}_3(l_1))^4  (\mathrm{sleaf}_3(l_2) )^2 +  4 (\mathrm{sleaf}_3(l_1))^2  (\mathrm{sleaf}_3(l_2) )^4 \label{I5}
\end{equation}
The following equations are obtained by differentiating with respect to variable$l_1$:
\begin{equation}
 \frac{ \partial p_1(l_1,l_2)}{ \partial l_1} = \frac{ \partial \mathrm{sleaf}_3(l_1) }{ \partial l_1} \frac{ \partial \mathrm{sleaf}_3(l_2) }{ \partial l_2} -3  \mathrm{sleaf}_3(l_2)  (\mathrm{sleaf}_3(l_1) )^5
 \label{I6}
\end{equation}
\begin{equation}
 \frac{ \partial p_2(l_1,l_2)}{ \partial l_1} = 3 (\mathrm{sleaf}_3(l_1) )^2 \mathrm{sleaf}_3(l_2) \frac{ \partial \mathrm{sleaf}_3(l_1) }{ \partial l_1} - (\mathrm{sleaf}_3(l_2) )^3 \frac{ \partial \mathrm{sleaf}_3(l_1) }{ \partial l_1}
 \label{I7}
\end{equation}
\begin{equation}
 \frac{ \partial p_3(l_1,l_2)}{ \partial l_1} = 16 (\mathrm{sleaf}_3(l_1) )^3 (\mathrm{sleaf}_3(l_2))^2 \frac{ \partial \mathrm{sleaf}_3(l_1) }{ \partial l_1} +8 \mathrm{sleaf}_3(l_1) (\mathrm{sleaf}_3(l_2) )^4 \frac{ \partial \mathrm{sleaf}_3(l_1) }{ \partial l_1}
 \label{I8}
\end{equation}
The following equations are obtained by differentiating with respect to variable $l_2$:
\begin{equation}
 \frac{ \partial p_1(l_1,l_2)}{ \partial l_2} = \frac{ \partial \mathrm{sleaf}_3(l_1) }{ \partial l_1} \frac{ \partial \mathrm{sleaf}_3(l_2) }{ \partial l_2} -3  \mathrm{sleaf}_3(l_1)  (\mathrm{sleaf}_3(l_2) )^5
 \label{I9}
\end{equation}
\begin{equation}
 \frac{ \partial p_2(l_1,l_2)}{ \partial l_2} = - 3 (\mathrm{sleaf}_3(l_2) )^2 \mathrm{sleaf}_3(l_1) \frac{ \partial \mathrm{sleaf}_3(l_2) }{ \partial l_2} + (\mathrm{sleaf}_3(l_1) )^3 \frac{ \partial \mathrm{sleaf}_3(l_2) }{ \partial l_2}
 \label{I10}
\end{equation}
\begin{equation}
 \frac{ \partial p_3(l_1,l_2)}{ \partial l_2} = 16 (\mathrm{sleaf}_3(l_2) )^3 (\mathrm{sleaf}_3(l_1))^2 \frac{ \partial \mathrm{sleaf}_3(l_2) }{ \partial l_2} +8 \mathrm{sleaf}_3(l_2) (\mathrm{sleaf}_3(l_1) )^4 \frac{ \partial \mathrm{sleaf}_3(l_2) }{ \partial l_2}
 \label{I11}
\end{equation}
Using Eq. (\ref{I1}), the following equations are obtained by differentiating with respect to variable $l_1$:
\begin{equation}
\begin{split}
& \frac{ \partial g(l_1,l_2)}{ \partial l_1} = \frac{ \Big( 2  p_1(l_1,l_2) \frac{ \partial p_1(l_1,l_2)   }{ \partial l_1} + 2  p_2(l_1,l_2) \frac{ \partial p_2(l_1,l_2)}{ \partial l_1} \Big)  p_3(l_1,l_2)  }{ 2g(l_1,l_2) p_3(l_1,l_2)^2 } \\
& - \frac{  \Big( p_1(l_1,l_2)^2 + p_2(l_1,l_2)^2 \Big) \frac{ \partial p_3(l_1,l_2)   }{ \partial l_1}   }{ 2g(l_1,l_2) p_3(l_1,l_2)^2 } 
\label{I12}
\end{split}
\end{equation}
Using Eqs. (\ref{I6}) $\sim$ (\ref{I8}), the numerator in the Eq. (\ref{I12}) is expanded as:
\begin{equation}
\begin{split}
& \Big( 2 p_1(l_1,l_2) 
\frac{ \partial p_1(l_1,l_2)   }{ \partial l_1} 
+ 2 p_2(l_1,l_2) \frac{ \partial p_2(l_1,l_2)}{ \partial l_1} \Big) p_3(l_1,l_2)   
\\
& - \Big( p_1(l_1,l_2)^2 +p_2(l_1,l_2)^2   \Big) \frac{ \partial p_3(l_1,l_2)}{ \partial l_1} 
\\
&=\Big( 2 \mathrm{sleaf}_3(l_1) -8(\mathrm{sleaf}_3(l_1))^5 (\mathrm{sleaf}_3(l_2))^2 -24(\mathrm{sleaf}_3(l_1))^3 (\mathrm{sleaf}_3(l_2))^4 \\
& -8 \mathrm{sleaf}_3(l_1) (\mathrm{sleaf}_3(l_2))^6 -16 (\mathrm{sleaf}_3(l_1))^5  (\mathrm{sleaf}_3(l_2))^8    \Big) \frac{ \partial \mathrm{sleaf}_3(l_1)}{ \partial l_1} \\
&+\Big( 2 \mathrm{sleaf}_3(l_2) -8(\mathrm{sleaf}_3(l_1))^5 (\mathrm{sleaf}_3(l_2))^2 -24(\mathrm{sleaf}_3(l_1))^4 (\mathrm{sleaf}_3(l_2))^3 \\
& -8 (\mathrm{sleaf}_3(l_1))^6 \mathrm{sleaf}_3(l_2) -16 (\mathrm{sleaf}_3(l_1))^8  (\mathrm{sleaf}_3(l_2))^5    \Big) \frac{ \partial \mathrm{sleaf}_3(l_2)}{ \partial l_2} \\
\label{I13}
\end{split}
\end{equation}
Usingthe Eq. (\ref{I1}), the following equation is obtained by differentiating with respect to the variable  $l_2$:
\begin{equation}
\begin{split}
& \frac{ \partial g(l_1,l_2)}{ \partial l_2} = \frac{ \Big( 2  p_1(l_1,l_2) \frac{ \partial p_1(l_1,l_2)   }{ \partial l_2} + 2  p_2(l_1,l_2) \frac{ \partial p_2(l_1,l_2)}{ \partial l_2} \Big)  p_3(l_1,l_2)  }{ 2g(l_1,l_2) p_3(l_1,l_2)^2 } \\
& - \frac{  \Big( p_1(l_1,l_2)^2 + p_2(l_1,l_2)^2 \Big) \frac{ \partial p_3(l_1,l_2)   }{ \partial l_2}   }{ 2g(l_1,l_2) p_3(l_1,l_2)^2 } 
\label{I14}
\end{split}
\end{equation}
Using Eqs. (\ref{I9}) $\sim$ (\ref{I11}), the numerator in the Eq.  (\ref{I14}) is expanded to obtain the following relation:
\begin{equation}
\frac{ \partial g(l_1,l_2)   }{ \partial l_1}= \frac{ \partial g(l_1,l_2)   }{ \partial l_2}
\label{I15}
\end{equation}
The following equation is derived from the Eq. (\ref{I15}) (see Appendix J).
\begin{equation}
g(l_1,l_2)=g(l_1+l_2,0)
\label{I16}
\end{equation}
Using the initial condition $ \mathrm{sleaf}_3(0)=0$ and $\frac{ \partial \mathrm{sleaf}_3(0)   }{ \partial l}=1$, the function $g(l_1+l_2,0)$ is obtained as follows:
\begin{equation}
\{ g(l_1+l_2,0) \}^2=\frac{ \{ p_1(l_1+l_2,0) \}^2 }{p_3(l_1+l_2,0)}+\frac{ \{ p_2(l_1+l_2,0) \}^2 }{p_3(l_1+l_2,0)}=( \mathrm{sleaf}_3(l_1+l_2) )^2  \label{I17}
\end{equation}
Using Eqs.  (\ref{I1}), (\ref{I16}) and (\ref{I17}), the following equation is obtained.
\begin{equation}
\begin{split}
&( \mathrm{sleaf}_3(l_1+l_2) )^2= \{ g(l_1+l_2,0) \}^2=\{ g(l_1,l_2) \}^2 \\
& =\frac{ \left\{
\mathrm{sleaf}_3(l_1)\frac{\mathrm{\partial sleaf}_3(l_2)}{\mathrm{\partial} l_2}
+\mathrm{sleaf}_3(l_2) \frac{\mathrm{\partial sleaf}_3(l_1) }{\mathrm{\partial} l_1}
\right\}^2
}
{1+4(\mathrm{sleaf}_3(l_1))^4(\mathrm{sleaf}_3(l_2))^2+4(\mathrm{sleaf}_3(l_1))^2(\mathrm{sleaf}_3(l_2))^4} \\ 
& +\frac{\left\{
(\mathrm{sleaf}_3(l_1))^3 \mathrm{sleaf}_3(l_2) - \mathrm{sleaf}_3(l_1)	(\mathrm{sleaf}_3(l_2))^3 
\right\}^2}
{1+4(\mathrm{sleaf}_3(l_1))^4(\mathrm{sleaf}_3(l_2))^2+4(\mathrm{sleaf}_3(l_1))^2(\mathrm{sleaf}_3(l_2))^4}
\label{I18}
\end{split}
\end{equation}

\appendix
\def\thesection{Appendix}
\section{J}
\label{Appendix J}
The necessary and sufficient condition to satisfy $g(l_1,l_2)=g(l_1+l_2,0)$ is that $\frac{ \partial g(l_1,l_2)   }{ \partial l_1}= \frac{ \partial g(l_1,l_2)   }{ \partial l_2}$ holds. Function $h(x,y)$ is defined as follows.
\begin{equation}
h(x,y)=g(x+y,x-y) \label{J1}
\end{equation}
\begin{equation}
l_1=x+y \label{J2}
\end{equation}
\begin{equation}
l_2=x-y \label{J3}
\end{equation}
By differentiating the Eq. (\ref{J1}) equation with respect to $y$, the following equation is obtained.
\begin{equation}
\begin{split}
& \frac{ \partial h(x,y) }{ \partial y}=\frac{ \partial g(x+y,x-y) }{ \partial l_1} \frac{ \partial l_1 }{ \partial y} + \frac{ \partial g(x+y,x-y) }{ \partial l_2} \frac{ \partial l_2 }{ \partial y} \\
& = \frac{ \partial g(x+y,x-y) }{ \partial l_1} - \frac{ \partial g(x+y,x-y) }{ \partial l_2}
\label{J4}
\end{split}
\end{equation}
Therefore, if the equation $\frac{ \partial g   }{ \partial l_1}= \frac{ \partial g   }{ \partial l_2}$ holds, the following equation holds.
\begin{equation}
\frac{ \partial h(x,y)   }{ \partial y}=0 \label{J5}
\end{equation}
Using the Eq. (\ref{J5}), we find that $h(x,y)$ is a function of $x$ and not of $y$. Therefore, the following equation holds for any constant $a$ and $b$:
\begin{equation}
h(x,a)=h(x,b) \label{J6}
\end{equation}
Here, we set the following equation:
\begin{equation}
x=b=\frac{ l_1+l_2  }{ 2 }=0 \label{J7}
\end{equation}
\begin{equation}
a=\frac{ l_1-l_2  }{ 2 }=0 \label{J8}
\end{equation}
The following equation is obtained by using the Eqs. (\ref{J4}), (\ref{J7}), and (\ref{J8}):
\begin{equation}
h(x,a)=h \Big( \frac{ l_1+l_2  }{ 2 }, \frac{ l_1-l_2  }{ 2 } \Big)
=g \Big( \frac{ l_1+l_2  }{ 2 } + \frac{ l_1-l_2  }{ 2 },  \frac{ l_1+l_2  }{ 2 } - \frac{ l_1-l_2  }{ 2 }  \Big)=g(l_1,l_2)   \label{J9}
\end{equation}
\begin{equation}
h(x,b)=h \Big( \frac{ l_1+l_2  }{ 2 }, \frac{ l_1+l_2  }{ 2 } \Big)
=g \Big( \frac{ l_1+l_2  }{ 2 } + \frac{ l_1+l_2  }{ 2 },  \frac{ l_1+l_2  }{ 2 } - \frac{ l_1+l_2  }{ 2 }  \Big)=g(l_1+l_2,0)   \label{J10}
\end{equation}
The following equation is obtained by using the Eqs. (\ref{J6}), (\ref{J9}) and (\ref{J10}).
\begin{equation}
g(l_1,l_2)=g(l_1+l_2,0) \label{J11}
\end{equation}
Conversely, if the Eq. (\ref{J11}) holds, the following relational expression can be obtained by using Eqs. (\ref{J1}) and (\ref{J11}).
\begin{equation}
h(x,y)=g(x+y,x-y)=g(2x,0) \label{J12}
\end{equation}
Eq. (\ref{J12}) is differentiated with respect to variable $y$ to obtain the following equation:
\begin{equation}
\frac{ \partial h(x,y)   }{ \partial y}=\frac{ \partial g(2x,0)   }{ \partial y}=0 \label{J13}
\end{equation}
Further, the following equation is obtained by using the Eq. (\ref{J4}).
\begin{equation}
\frac{ \partial g(x+y,x-y) }{ \partial l_1} = \frac{ \partial g(x+y,x-y) }{ \partial l_2}
\label{J14}
\end{equation}
Using the Eqs. (\ref{J2}) and (\ref{J3}), the following equation is obtained.
\begin{equation}
\frac{ \partial g(l_1,l_2) }{ \partial l_1} = \frac{ \partial g(l_1,l_2) }{ \partial l_2}
\label{J15}
\end{equation}

\nocite{*}

\bibliographystyle{unsrt}  
\bibliography{ref}  %%% Remove comment to use the external .bib file (using bibtex).

\begin{thebibliography}{10}

\bibitem{Kaz_cl}
Kazunori Shinohara.
\newblock Special function: Leaf function $r=cleaf_n(l)$ (second report).
\newblock {\em Bulletin of Daido University}, 51:39--68, 2015.

\bibitem{Kaz_sl}
Kazunori Shinohara.
\newblock Special function: Leaf function $r=sleaf_n(l)$ (first report).
\newblock {\em Bulletin of Daido University}, 51:23--38, 2015.

\bibitem{Gauss}
C.F. Gauss, W.C. Waterhouse, and A.A. Clarke.
\newblock {\em Disquisitiones Arithmeticae}.
\newblock Springer-Verlag, 1986.

\bibitem{Car}
B.~C. Carlson.
\newblock Algorithms involving arithmetic and geometric means.
\newblock {\em The American Mathematical Monthly}, 78(5):496--505, 1971.

\bibitem{Neu}
Edward Neuman.
\newblock On lemniscate functions.
\newblock {\em Integral Transforms and Special Functions}, 24(3):164--171,
  2013.

\bibitem{Wei}
E.W. Weisstein.
\newblock {\em CRC Concise Encyclopedia of Mathematics, Second Edition}.
\newblock CRC Press, 2002.

\bibitem{Roy}
R.~Roy.
\newblock {\em Elliptic and Modular Functions from Gauss to Dedekind to Hecke}.
\newblock Cambridge University Press, 2017.

\bibitem{Ayo}
Raymond Ayoub.
\newblock The lemniscate and fagnano's contributions to elliptic integrals.
\newblock {\em Archive for History of Exact Sciences}, 29(2):131--149, 1984.

\bibitem{Euler}
L.~Euler.
\newblock {\em Leonhardi Euleri Opera omnia}.
\newblock Teubner, Leipzig, 1911.

\bibitem{Jacobi}
C.~G.~J. Jacobi.
\newblock {\em Opuscula Mathematica. Mathematische Werke}.
\newblock Nabu Press, 2010.

\bibitem{Kaz_add}
Kazunori Shinohara.
\newblock Addition formulas of leaf functions according to integral root of
  polynomial based on analogies of inverse trigonometric functions and inverse
  lemniscate functions.
\newblock {\em Applied Mathematical Sciences}, 11:2561--2577, 01 2017.

\bibitem{Sibanda}
Precious Sibanda and Ahmed Khidir.
\newblock A new modification of the hpm for the duffing equation with cubic
  nonlinearity.
\newblock {\em International Conference on Applied and Computational
  Mathematics}, 2, 01 2011.

\bibitem{Dai}
Honghua Dai.
\newblock A simple collocation scheme for obtaining the periodic solutions of
  the duf?ng equation, and its equivalence to the high dimensional harmonic
  balance method: Subharmonic oscillations.
\newblock {\em Computer Modeling in Engineering and Sciences}, 84:459--497, 09
  2012.

\bibitem{Bulbul}
Berna Bulbul.
\newblock Numerical solution of duffing equation by using an improved taylor
  matrix method.
\newblock {\em Journal of Applied Mathematics}, 2013, 06 2013.

\bibitem{ELIASZUNIGA2014849}
Alex Elias-Zuniga.
\newblock Quintication method to obtain approximate analytical solutions of
  non-linear oscillators.
\newblock {\em Applied Mathematics and Computation}, 243:849 -- 855, 2014.

\bibitem{Sayevand}
Khosro Sayevand, Dumitru Baleanu, and M.~Fardi.
\newblock A perturbative analysis of nonlinear cubic-quintic duffing
  oscillators.
\newblock {\em Proceedings of the Romanian Academy Series A - Mathematics
  Physics Technical Sciences Information Science}, 15:228--234, 07 2014.

\bibitem{Novin}
Reza Novin and Ziba Dastjerd.
\newblock Solving duffing equation using an improved semi-analytical method.
\newblock {\em Communications on Advanced Computational Science with
  Applications}, 2015:54--58, 01 2015.

\bibitem{Zhang}
Ying Zhang, Lin Du, Xiaole Yue, Qun Han, and Tong Fang.
\newblock Analysis of symmetry breaking bifurcation in duffing system with
  random parameter.
\newblock {\em CMES - Computer Modeling in Engineering and Sciences},
  106:37--51, 05 2015.

\bibitem{Stoyanov}
Svetlin Stoyanov.
\newblock Analytical and numerical investigation on the duffing oscilator
  subjected to a polyharmonic force excitation.
\newblock {\em Journal of Theoretical and Applied Mechanics}, 45, 03 2015.

\bibitem{ELNAGGAR20161581}
A.M. El-Naggar and G.M. Ismail.
\newblock Analytical solution of strongly nonlinear duffing oscillators.
\newblock {\em Alexandria Engineering Journal}, 55(2):1581 -- 1585, 2016.

\bibitem{Weli}
Majeed Weli and Sayl Abd-Al-Razaq.
\newblock A semi analytical iterative technique for solving duffing equations.
\newblock {\em International Journal of Pure and Applied Mathematics},
  108:871--885, 08 2016.

\bibitem{Hosen}
Md~Hosen.
\newblock Solution of cubic-quintic duffing oscillators using harmonic balance
  method.
\newblock {\em Malaysian Journal of Mathematical Sciences}, 10:181--192, 02
  2016.

\bibitem{CHOWDHURY20173962}
M.S.H. Chowdhury, Md.~Alal Hosen, Kartini Ahmad, M.Y. Ali, and A.F. Ismail.
\newblock High-order approximate solutions of strongly nonlinear cubic-quintic
  duffing oscillator based on the harmonic balance method.
\newblock {\em Results in Physics}, 7:3962 -- 3967, 2017.

\bibitem{Karahan}
Mustafa Karahan and Mehmet Pakdemirli.
\newblock Free and forced vibrations of the strongly nonlinear cubic-quintic
  duffing oscillators.
\newblock {\em Zeitschrift fur Naturforschung A}, 72, 01 2017.

\bibitem{KOVACIC20161}
Ivana Kovacic, Livija Cveticanin, Miodrag Zukovic, and Zvonko Rakaric.
\newblock Jacobi elliptic functions: A review of nonlinear oscillatory
  application problems.
\newblock {\em Journal of Sound and Vibration}, 380:1 -- 36, 2016.

\bibitem{ELIASZUNIGA20132574}
Alex Elias-Zuniga.
\newblock Exact solution of the cubic-quintic duffing oscillator.
\newblock {\em Applied Mathematical Modelling}, 37(4):2574 -- 2579, 2013.

\bibitem{Belendez}
Augusto Belendez, T.~Belendez, Francisco Martinez~Guardiola, Carolina
  Villalobos, M.~Alvarez, and Enrique Arribas.
\newblock Exact solution for the unforced duffing oscillator with cubic and
  quintic nonlinearities.
\newblock {\em Nonlinear Dynamics}, 86, 08 2016.

\bibitem{Nwamba}
Joshua Nwamba.
\newblock Application of iteration perturbation method and variational
  iteration method to a restrained cargo system modeled by cubic-quintic-septic
  duffing equation.
\newblock {\em International Journal of Mechanics and Applications},
  2013:63--69, 03 2013.

\bibitem{Serge}
Serge~Bruno Yamgoue, Olivier~Tiokeng Lekeufack, and Timoleon~Crepin Kofane.
\newblock Harmonic balance for non-periodic hyperbolic solutions of nonlinear
  ordinary differential equations.
\newblock {\em Mathematical Modelling and Analysis}, 22(2):140--156, 2017.

\bibitem{REMMI20182085}
S.K. Remmi and M.M. Latha.
\newblock Cubic quintic septic duffing oscillator: An analytical study.
\newblock {\em Chinese Journal of Physics}, 56(5):2085 -- 2094, 2018.

\bibitem{Koudahoun}
Herve Koudahoun, Y.~Kpomahou, Judith Akande, and Kolawole Kegnide~Damien Adjai.
\newblock Chaotic dynamics of an extended duffing oscillator under periodic
  excitation.
\newblock {\em World Journal of Applied Physics}, 3:34--50, 08 2018.

\bibitem{Kaz_ch}
Kazunori Shinohara.
\newblock Special function: Hyperbolic leaf function $r=cleafh_n(l)$ (second
  report).
\newblock {\em Bulletin of Daido University}, 52:83--105, 2016.

\bibitem{Kaz_sh}
Kazunori Shinohara.
\newblock Special function: Hyperbolic leaf function $r=sleafh_n(l)$ (first
  report).
\newblock {\em Bulletin of Daido University}, 52:65--81, 2016.

\end{thebibliography}

\end{document}